\documentclass[a4paper,12pt]{article}

\makeatletter
\let\cl@chapter\undefined
\makeatletter

\usepackage{mathrsfs,amssymb,amsmath,amsthm}
\usepackage{authblk}

\usepackage{graphicx}
\usepackage[dvipsnames]{xcolor}
\usepackage{anyfontsize}
\usepackage{tikz}
\usetikzlibrary{shapes,positioning}
\usepackage[top=40pt, bottom=60pt,left=70pt, right=70pt]{geometry}

\usepackage{soul}

\usepackage[pagebackref,colorlinks,citecolor=Plum,urlcolor=Periwinkle,linkcolor=DarkOrchid]{hyperref}
\usepackage[square,numbers]{natbib}
\usepackage{verbatim} 
\RequirePackage[labelsep=space,tableposition=top]{caption}

\usepackage{subcaption}

\usepackage{float}
\usepackage{subcaption}
\usepackage{cleveref} 
\usepackage{enumitem}
\restylefloat{figure}

\usepackage{booktabs} 
\usepackage{multirow}

\newcommand{\pv}[2]{#1 a^\perp \otimes a^\perp + #2 n^\perp \otimes n^\perp}

\newcommand{\re}{\ensuremath{\mathbb{R}}}
\newcommand{\ten}[2]{#1_{#2} \otimes #1_{#2}}
\newcommand{\m}{\re^{2\times 2}}
\newcommand{\msym}{\re^{2\times 2}_{sym}}

\newcommand{\pro}[2]{\left \langle #1 \, , #2 \right \rangle}
\DeclareMathOperator{\rpartial}{\partial_{ri}}
\DeclareMathOperator{\Aff}{Affine}
\DeclareMathOperator{\rank}{Rank}
\DeclareMathOperator{\Id}{Id}

\DeclareMathOperator{\relint}{rel\, int}
\DeclareMathOperator{\diag}{diag}

\newtheorem{theorem}{Theorem}[section]
\newtheorem{proposition}{Proposition}[section]

\newtheorem{remark}{Remark}[section]
\newtheorem{lemma}{Lemma}[section]

\newtheorem*{theorem*}{Theorem}
\newtheorem*{proposition*}{Proposition}
\newtheorem*{remark*}{Remark}
\newtheorem*{lemma*}{Lemma}

\everymath{\displaystyle}

\begin{document}

\date{January 2021}

\author[1]{Capella, Antonio}
\author[2]{Morales, Lauro}
\affil[1]{\small Instituto de Matem\'aticas-UNAM \\ lmm@ciencias.unam.mx
}
\affil[2]{\small Instituto de Matem\'aticas-UNAM\\
capella@im.unam.mx}

\title{On the symmetric lamination convex and quasiconvex hull for the coplanar $n$-well problem in two dimensions.}
\maketitle
\begin{abstract}
\hskip -.2in
\noindent
We study some particular cases of the $n$-well problem in two-dimensional linear elasticity.  Assuming that every well in $\mathcal{U}\subset\mathbb{R}^{2\times 2}_\text{sym}$ belong to the same two-dimensional affine subspace, we characterize the symmetric lamination convex hull $L^e(\mathcal{U})$
for any number of wells in terms of the symmetric lamination convex hull of all three-well subsets contained in $\mathcal{U}$.
For a family of four-well sets where two pairs of wells are rank-one compatible, we show that the symmetric lamination convex and quasiconvex hulls coincide, but are strictly contained in its convex hull $C(\mathcal{U})$.  We extend this result to some particular configurations of $n$ wells.  Most of the proofs are constructive, and we also present explicit examples. 
\end{abstract}

\section{Introduction}\label{intro}
In this paper, we consider a nonconvex functional of the form
\begin{equation}\label{eq:functional}
I(u)=\int_\Omega \phi(e(D u)) dx, 
\end{equation}
where $u:\Omega\to \re^2$ is the displacement, $\Omega$ is a bounded region in $\re^2$, $e(D u):=(D u + (D u)^T)/2$ denotes its symmetrized deformation gradient or linear strain, and $\phi$ has a multi-well structure. 
Namely,  $\phi(M+S)=\phi(M)$ for every $M\in \msym$, and $S\in \re^{2\times 2}_{skew}$,  
$\phi(U)=0$ for every $U$ in the finite set of wells $\mathcal{U} = \{U_1,\, U_2,\, \cdots U_n\}\subset \msym$, and $\phi>0$ otherwise. 
This type of functionals typically appears in problems of hyper-elastic material in the framework of linear elasticity theory, see \cite{K,TZ2008}.

As in  \cite{BFJK,Z2004}, we say that two different matrices $M_1,\,M_2\in \msym$  are  {\em compatible}  if there exists a skew--symmetric matrix $S$ such that   $\rank(M_1-M_2+S) \le 1$, otherwise $M_1$ and $M_2$ are {\em incompatible}.
We also say that $M_1$ and $M_2$ are {\em rank-one compatible} provided $\rank (M_1-M_2)\leq 1$, or equivalently, whenever $\det (M_1-M_2)=0$.
It readily follows that if $M_1$ and $M_2$ are compatible, then $M_1-M_2+S=a\otimes n$  for some $a\in \re^2$ and $n\in \mathbb{S}^1$. 
Hence, if in addition $M_1,\,M_2\in \mathcal{U}$, then layered structures, or simple laminates, can be constructed and nontrivial minimizing sequences  of \cref{eq:functional} emerge in suitable functional spaces, see \cite{BFJK,RJ}. 
In this work, we are interested in characterizing the effective linear strains generated by microstructures, namely minimizing sequences.
For this characterization, we need some notions of semi-convex analysis that we recall in the next paragraphs.

We assume standard definitions of semiconvex functions, see \cite{BDa}. Let $n_m=C(2m,m)-1,$ and $C(2m,m)$ be the binomial coefficient between $2m$ and $m$. A function $f:\re^{m\times m} \rightarrow \re$ is {\em polyconvex} 
if there exists a convex function $G:\re^{m_n}\rightarrow \re$, such that 
$$f(M) = G\circ T(M),$$
where $T:\re^{m\times m}\rightarrow \re^{n_m}$  is given by  
$M\mapsto (M, adj_2(M),\ adj_3(M),\ \cdots, \det M ).$
The notation $adj_k(M)$ stands for the matrix of all $k\times k$ minors of $M$.  Also,  $f$ is {\em rank-one convex} if for each pair of rank-one compatible matrices $M_1,M_2\in \re^{m\times m}$, 
\[
f(\lambda M_1 +(1-\lambda)M_2) \leq \lambda f(M_1) +(1-\lambda)f(M_2), \ \ \lambda\in [0,1].
\]
Finally, $f$ is {\em quasiconvex} provided it is a Borel measurable and locally bounded function such that
\[
f(M) \leq \inf_{\substack{\phi \in C^{\infty}_0(\Omega,\re^m)}} \ \frac{1}{|\Omega|}\int_{\Omega} f(M + D \phi) \,dx.
\]
The function $f:\re^{m\times m}_{sym} \rightarrow \re$ is \emph{symmetric semi-convex} if $f(e(\cdot)):\re^{m\times m} \rightarrow \re$ is semi-convex.
As stated in \cite{Anja19,Z2002}, it follows that $f:\re^{m\times m}_{sym} \rightarrow \re$ is  {\em symmetric quasiconvex} if and only of for every $U\in \re^{m\times m}_{sym}$ 
\[
f(U)\leq \inf\left\{\frac{1}{|\Omega|} \int_\Omega f(U +e(D\phi))dx \, \middle | \, \phi\in C^\infty_0(\Omega, \re^{m} ) \right\}.
\] 
Also $f$  is  {\em symmetric rank-one convex} if for every two compatible matrices $U_1$ and $U_2\in \re^{m\times m}_{sym}$ and $\lambda\in [0,1]$,
\[
f(\lambda U_1 +(1-\lambda)U_2) \leq \lambda f(U_1) +(1-\lambda)f(U_2). 
\] 

For any compact set $\mathcal{U}\in \re^{n\times n}_{sym}$, its \emph{symmetric semi-convex hulls} are defined by mean of cosets \cite{Z2002}, more precisely,
\begin{equation}\label{semihulls}
    S^e(\mathcal{U}) = \left\{A\in \re^{n\times n}\, | \, f(A) \leq \sup_{B\in \mathcal{U}} f(B), \ f \ \mbox{symmetric semi-convex}\right\},
\end{equation}
where $S$ must be replace by $R,Q$ and $P$ for the symmetric rank-one convex, quasiconvex and polyconvex hull of the set $\mathcal{U}$. 
The set $Q^e(\mathcal{U})$ is relevant \cite{Z2002} since it corresponds to the set of effective linear strains generated by the microstructures with symmetric deformation gradient in $\mathcal{U}$.
As in the nonlinear case, to determine the symmetric quasiconvex hull of a compact set $\mathcal{U}$ is a challenging task \cite{Lazlo2006,Z2002}. 
Even for particular choices of $\mathcal{U}$ \cite{HK2017}, the  rank-one convex and quasiconvex hulls are difficult to compute. 
Explicit examples of quasiconvex hulls \cite{BD,BFJK,Bt,Z1998} are scarce, and  most of the nontrivial examples are not known explicitly. 

An inner bound for the symmetric quasiconvex hull is given by $L^e(\mathcal{U})$, the symmetric lamination convex hull of $\mathcal{U}$, see Zhang \cite{A2016,Z2002}. This set is the union of all compatible lamination hulls of every ranks, namely
\[
L^e(\mathcal{U}) = \bigcup_{i=0}^\infty L^{e,i}(\mathcal{U}),
\]
where $L^{e,0}(\mathcal{U}) = \mathcal{U},$ and 
$$
L^{e,i+1}(\mathcal{U}) = \{\lambda A + (1-\lambda) B \in\re^{n\times n}_{sym} \, | \, \lambda \in [0,1]\mbox{ and } A,B\in L^{e,i}(\mathcal{U}) \mbox{ are compatible}\}. 
$$
It is also known \cite{Anja19,Z2002} that for any compact set $\mathcal{U}\subset \re^{n\times n}_{sym}$,
\begin{equation}
\mathcal{U} \subset L^{e}(\mathcal{U}) \subset R^e(\mathcal{U}) \subset Q^e(\mathcal{U}) \subset P^e(\mathcal{U}) \subset C(\mathcal{U}). 
\end{equation}

In the present paper, we are interested in the particular case where the set of wells belong to an affine subspace $\Pi_Q\subset \msym$ with codimension one, briefly called a coplanar set. To be more precise,  we say that $\{U_1,U_2,\,\cdots, U_n\}\subset \re^{2\times 2}_{sym}$ is 
coplanar, if there exists $Q\in \re^{2\times 2}_{sym}$ and $\delta\in \re$ 
fix such that $\pro{Q}{U_i} = \delta$ for every $i=1,2,\cdots,n$. Here,  $<\cdot,\cdot>:\re^{2\times 2}_{sym}\times \re^{2\times 2}_{sym}\rightarrow \re$  denotes the Frobenius inner product, and we call $Q$ the normal to the affine space $\Pi_Q$.


In \cite{camo2020}, the authors provide a characterization of the lamination convex hull of the three-well case in two-dimensions depending on the wells' compatibility. 
As part of the paper's arguments, the authors characterized the set of incompatible wells with a given well $U$. More precisely, for every $U\in\msym$, 
the set 
\begin{equation}\label{eq:cone}
\mathcal{C}_U := \left\{V\in\msym \mbox{ such that } \  \|V-U\| < |\pro{V-U}{\Id}| \right \}.
\end{equation}
is the incompatible cone of $U$.  That is, $U$ and $V$ are incompatible if and only if $V\in\mathcal{C}_U$. 
The axis of the cone $\mathcal{C}_U$ is aligned to the subspace $\{t\,\Id \, |\, t\in \re\}$.
Moreover, the sets' boundary $\partial \mathcal{C}_0$ consists of all rank-one symmetric matrices (see lemma~2 in \cite{camo2020} and \cref{fig:cones}  for a representation under the isomorphism between $\re^3$ and  $\msym$).
We also showed in \cite{camo2020} that, if $\mathcal{U}\subset\msym$ is a set of three linearly independent matrices and there exist at least two of them that are rank-one compatible, then $Q^e(\mathcal{U}) = L^e(\mathcal{U})$. Moreover, only in the case where  all wells are compatible $L^e(\mathcal{U})=C(\mathcal{U})$, see also \cite{Bt,A2016}. In this paper, we extend this type of results for more than three wells.

\section{Main results}

 We base part of our present paper's results on the three-well characterization of the symmetric lamination hull from \cite{camo2020}.  
Let $\mathcal{U}$ be a finite coplanar set of wells,  we define  $\mathcal{F}$ as the family of triplets contained in $\mathcal{U}$, namely  
\begin{equation}\label{laminarcontation}
\mathcal{F}=\{\{U_i,U_j,U_k\}\,| \, i,j,k\in \{1,2,\dots,n\} \mbox{ are diferent indexes} \}.
\end{equation}
Furthermore, since $\mathcal{V}\subset \mathcal{U}$ for every $\mathcal{V}\in \mathcal{F}$ it follows that 
$L^{e,i}(\mathcal{V}) \subset L^{e,i}(\mathcal{U})$ for each $i\in \mathbb{N}$, and 
\begin{equation}\label{laminarcontation2}
\bigcup_{\mathcal{V}\in \mathcal{F}}L^{e,i}(\mathcal{V}) \subset L^{e,i}(\mathcal{U}).
\end{equation}
We characterize the set $L^e(\mathcal{U})$, with $\mathcal{U}$ a set of $n$ coplanar  wells in our first result.
\begin{theorem}[Laminar convex]\label{thm:lamconv}
Let $\mathcal{U}\subset \msym$, be a finite coplanar set and let $\mathcal{F}$ be as in \cref{laminarcontation}, then
\[
L^e(\mathcal{U})=\bigcup_{\mathcal{V}\in \mathcal{F}}L^e(\mathcal{V}).
\]
\end{theorem}
Theorem~\ref{thm:lamconv} combined with the characterization of \cite{camo2020} for three-wells (see also Proposition~\ref{P:lamconv} below) give a full description of the symmetric laminar convex hull for any number of coplanar wells. 

\begin{figure}[ht]
\begin{centering}
\includegraphics[scale=0.5]{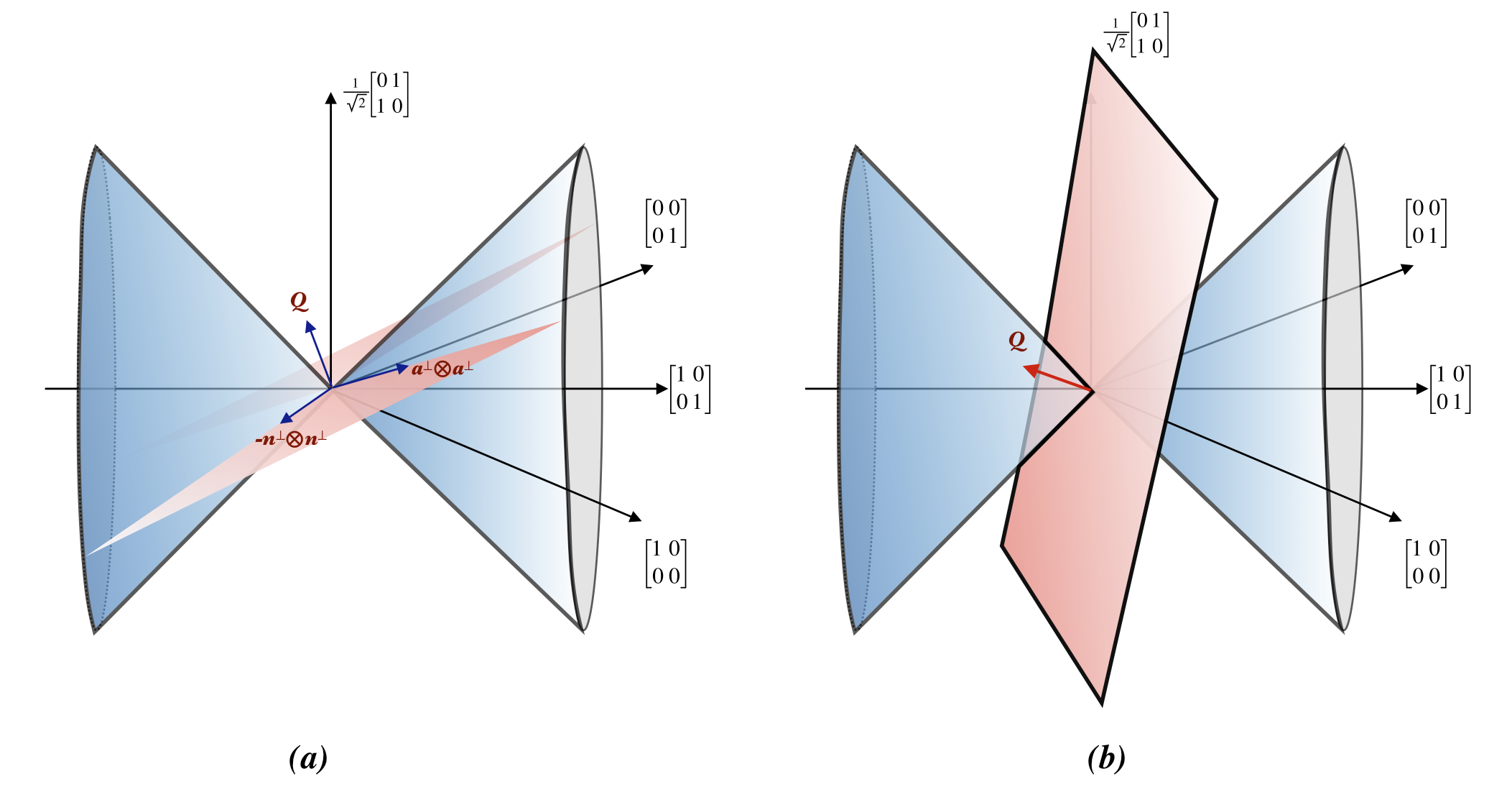}
\caption{Figure (a) shows the plane $\Pi_Q$ with $\det Q < 0$ and the incompatible cone center at the origin. 
The zero matrix and $U\in \Pi_Q$ are compatible if and only if $U$ belongs to the red regions. Figure (b) displays the case where $\det Q>0$ and every point in the plane is compatible with the origin.}
\label{fig:cones}
\end{centering}
\end{figure}

In \cite{camo2020}, we showed that if the normal to $\Pi_Q$ is such that $\det Q\geq 0$, then $\Pi_Q$ does not intersect the incompatible cone of the elements in $\mathcal{U}$ outside the vertices, and every pair of wells in $\mathcal{U}$ are pairwise compatible. Thus, as already mentioned, in general, $L^e(\mathcal{U}) = Q^e(\mathcal{U})=C(\mathcal{U})$ \cite{Bt,A2016}. Therefore, the interesting case is when $\det Q < 0$, and $\Pi_Q$ intersects the incompatible cone $\mathcal{C}_{U}$ centered at $U$, for every element in $\mathcal{U}$ (see Figure~\ref{fig:cones}).

We introduce some notations before presenting our next result. First,  for any coplanar set $\mathcal{U}$ we denote by $\mathscr{C}_Q(U)={\mathcal C}_U^c\cap \Pi_Q$ the planar \emph{compatible} cone of $U$.
In \cite{camo2020} (see also \cref{planeQV} in \cref{S:section3}), we show that if $\det Q<0$, then there exist $a,\, n\in S^1$, linearly independent such that $\pro{Q}{a^\perp \otimes a^\perp} = \pro{Q}{n^\perp \otimes n^\perp}=0$, and  $\{a^\perp \otimes a^\perp, n^\perp \otimes n^\perp\}$ is a  basis for the  subspace  $\Pi_Q$ (see Figure~\ref{fig:cones} (a)). 
Here, $v^\perp$  stands for the $\pi/2$ counter-clockwise rotation of any vector $v\in \re^2$.
With this in hand, we define the {\em upper} and {\em lower parts} of $\mathscr{C}_Q(U)$ as 
\begin{equation*}
\mathscr{C}^+_Q(U)=\{ V\in\re^{2\times 2}_{sym} \ \text{ such that }\  
V= U+\xi a^\perp \otimes a^\perp + \eta n^\perp \otimes n^\perp, \mbox{ for some }  \xi\geq 0 \geq \eta   \}
\end{equation*}
and
\begin{equation*}
\mathscr{C}^-_Q(U)=\{ U\in\re^{2\times 2}_{sym} \ \text{ such that }\  
V= U+ \xi a^\perp \otimes a^\perp + \eta n^\perp \otimes n^\perp, \mbox{ for some }  \eta\geq 0 \geq \xi   \},    
\end{equation*}
respectively. 



\begin{figure}[ht]
\begin{centering}
\includegraphics[scale=0.3]{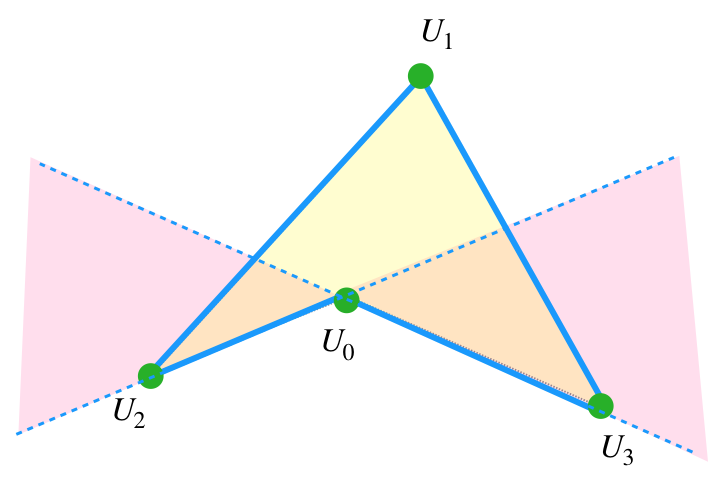}
\caption{Four coplanar wells in a wedge configuration}
\label{fig:wedge}
\end{centering}
\end{figure}

Second, we identify a particular geometric configuration of the set $\mathcal{U}$ in the case of four wells that we need to rule out in our results. 
We say that the set $\mathcal{U} = \{U_0,U_1,U_2,U_3\}\subset \re^{2\times 2}_{sym}$ is in a {\em wedge configuration} if there exists a subset of three wells, say $\mathcal{V}= \{U_1,U_2,U_3\}$,  such that 
there is only one incompatible pair of wells in $\mathcal{V}$, and 
the remaining well $U_0\in \relint C(\mathcal{V})$ is rank-one compatible with each element in the incompatible pair of wells in $\mathcal{V}$, see \cref{fig:wedge}.

Now we consider a set $\mathcal{U}$ of four coplanar wells with $\det Q <0$. 
In this case, we prove that the laminar convex hull of $\mathcal{U}$ equals $Q^e(\mathcal{U})$ for a family of four wells with two pairs of rank-one compatible wells. Moreover, these symmetric semi-convex hull are strictly contained in $C(\mathcal{U})$.

\begin{theorem}[Four coplanar wells]\label{thm:4wells}
Let $\mathcal{U}\subset \re^{2\times 2}_{sym}$ be a set of four coplanar wells such that all its elements have at least another compatible one, its plane's normal satisfies $\det Q <0$, and $\mathcal{U}$ is not in a wedge configuration. Furthermore, assume there exist two different subsets $\{V_1,V_2\},\{W_1,W_2\}\subset \mathcal{U}$ of rank-one compatible pairs and let $D = C(\{W_1,W_2\}) \cup C(\{V_1,V_2\})$. 
If any of the following conditions holds,
\begin{enumerate}
\item\label{it:le1} The set $D$ is disconnected,
\item\label{it:le3} $D$ is a connected set and $\mathcal{U}\subset D$,
\item\label{it:le2} The intersection of the sets $\{V_1,V_2\}$ and $\{W_1,W_2\}$ has only one element, say $V$, and $D$ is contained either in the upper or in the lower part of $\partial(\mathscr{C}_Q(V))$,
\end{enumerate}
then  $L^e(\mathcal{U}) = Q^e(\mathcal{U})$. 
\end{theorem}

Theorem~\ref{thm:4wells} extends the results for 
the three-well case presented in \cite{camo2020}. 
It is not a full characterization of the quasiconvex hull for the four-well problem, but in combination with Theorem~\ref{thm:lamconv} it let us produce examples of set with four and more wells, where we explicitly compute $L^e(\mathcal{U})$ and $Q^e(\mathcal{U})$.  
Besides the wedge configurations, we have only two cases that are not included in the statement 
of \cref{thm:4wells}. In particular, when the intersection of the sets $\{V_1,V_2\}$ and $\{W_1,W_2\}$ has only one element $V$, and $D$ is aligned to the left and right parts of $\mathscr{C}_Q(V)$.

\begin{figure}[ht]
     \centering
     \begin{subfigure}[b]{0.6\textwidth}
         \centering
\includegraphics[height=0.15\textheight]{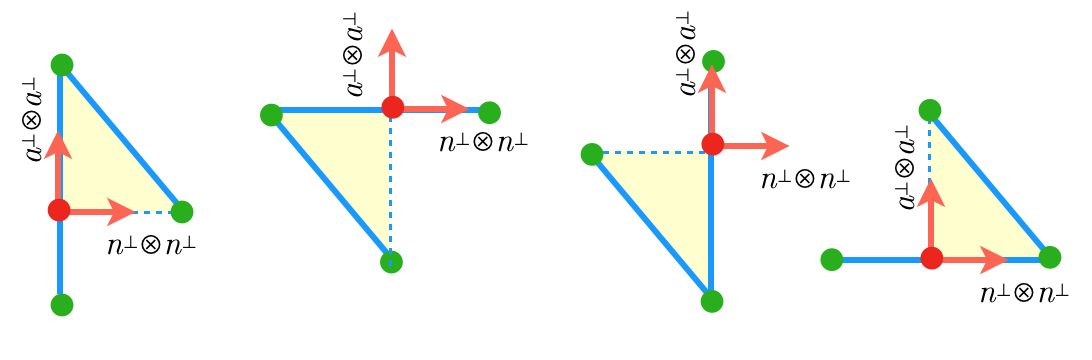}
         \caption{Basic three-well blocks}
         \label{fig:bb3}
     \end{subfigure}
     \hspace{0.8cm}
     \begin{subfigure}[b]{0.3\textwidth}
         \centering
         \includegraphics[height=0.15\textheight]{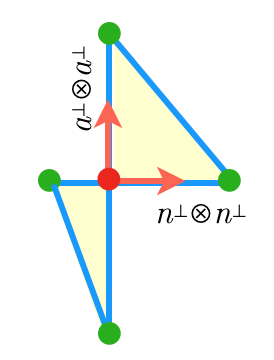}
         \caption{Basic four-well block}
         \label{fig:bb4}
     \end{subfigure}
     \caption{The five basic block configurations. Green dots represent the wells in $\mathcal{U}$ while 
     red dots represent the center of each basic block.
     To keep figures simple we assume that the angle between $a^\perp\otimes a^\perp$ and $n^\perp\otimes n^\perp$ directions is $\pi/2$.}
\end{figure}

In our final main result, we give a constructive procedure for sets $\mathcal{U}$ with an arbitrary finite number of coplanar wells where $L^e(\mathcal{U})=Q^e(\mathcal{U})$, but they are strictly contained in $C(\mathcal{U})$.
 We set up these configurations in terms of five \emph{basic blocks} with three or four wells.   
 Three-well basic blocks have one compatible, one rank-one compatible, and one incompatible pair among its wells,  see \cref{fig:bb3}, we also refer to this type of basic block as {\em flag configurations}.  The four-well basic block corresponds to the case of \cref{it:le3} in \cref{thm:4wells}, see \cref{fig:bb4}. Notice that in each basic block, the segment generated by a  rank-one pair of wells contains a well $U_0$ that is rank-one compatible with any other element in the flag configuration. We called $U_0$ the {\em basic block's center}. By \cref{thm:3well} and \cref{thm:4wells}, if $\mathcal{U}$ is a basic block, then $L^e(\mathcal{U})=Q^e(\mathcal{U})$.

\begin{figure}[ht]
     \centering
     \begin{subfigure}[b]{0.4\textwidth}
         \centering
\includegraphics[height=0.1\textheight]{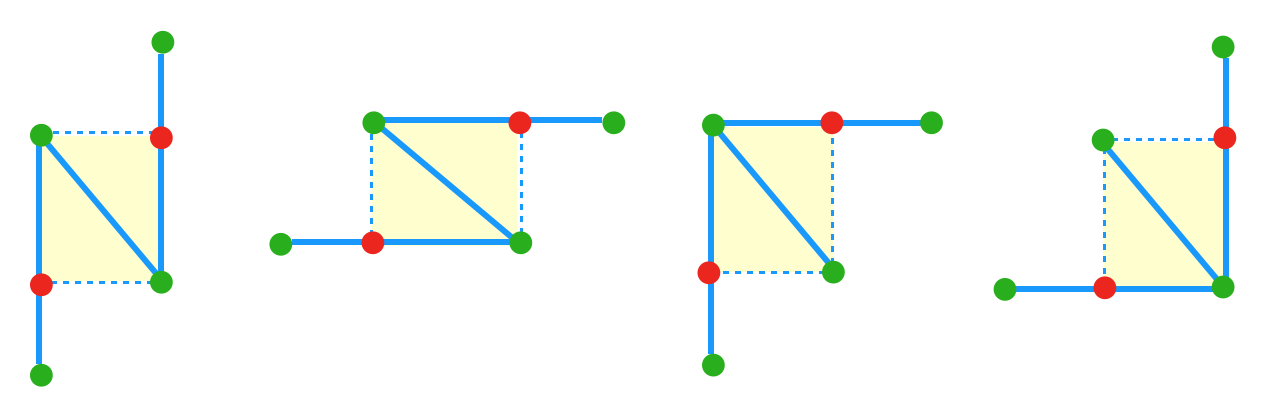}
         \caption{Four-well basic configuration}
         \label{fig:bc2}
     \end{subfigure}
     \hspace{1.5cm}
     \begin{subfigure}[b]{0.4\textwidth}
         \centering
         \includegraphics[height=0.1\textheight]{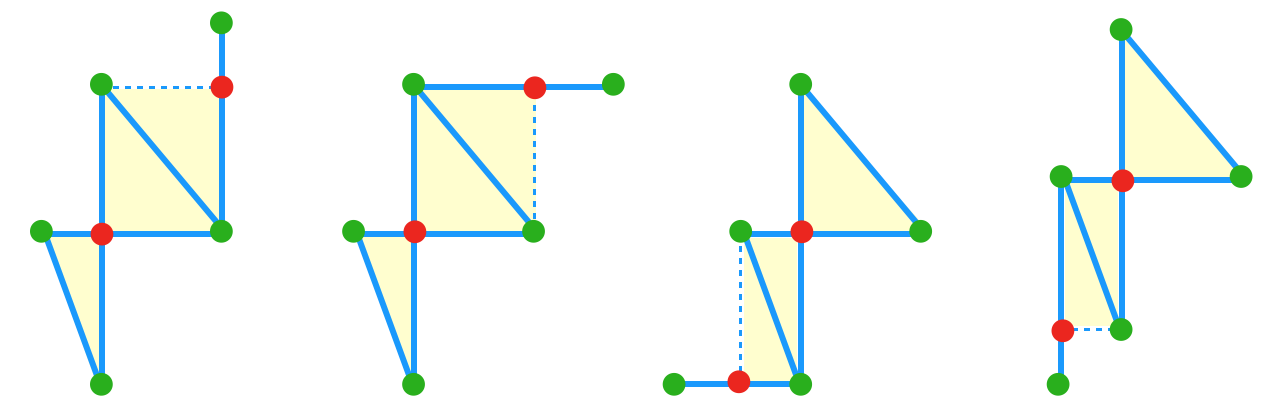}
         \caption{Five-well basic configuration}
         \label{fig:bc2-1}
     \end{subfigure}
     \caption{Basic configurations obtained by sticking two basic blocks. The configurations presented in \cref{fig:bc2} are included in the conditions of \cref{thm:4wells}.}
     \label{fig:bc2all}
\end{figure}

The basic blocks can be stick together to get other configurations with more than three wells where the same affirmation hold.    
We say $\mathcal{U}$ is a  \emph{basic configuration} if there are no more than two colinear wells in the set and it is the union of finitely many adjacent basic blocks, where every pair of adjacent basic blocks share two compatible (but not rank-one compatible) wells. In \cref{fig:bc2all} we show some examples of basic configurations.  Notice that, if $\mathcal{U}$ is a basic configuration, then $\mathcal{U}=\mathcal{U}_1\cup \cdots \cup \mathcal{U}_n$, where $\mathcal{U}_i$ is a basic block for every $i=1,2,\cdots, n$; and if $i\neq n$, then the set $\mathcal{U}_i\cap\mathcal{U}_{i+1}=\{V_i,W_i\}$ where $\det(V_i-W_i)<0$. 
By \cref{prop:bc} (see below), there are at most two three-wells basic blocks in any basic configuration. 
Furthermore, these three-well basic blocks appear only as the first one or the last one in $\mathcal{U}$.

\begin{theorem}[Basic configurations]\label{thm:iterative}
Let $\mathcal{U}\subset \re^{2\times 2}_{sym}$ be a finite coplanar set of  wells. If $\mathcal{U}$ is a basic configuration, then $L^e(\mathcal{U})=Q^e(\mathcal{U})$.
\end{theorem}

In the proof of \cref{thm:iterative} we assume that many pairs of wells in $\mathcal{U}$ are rank-one compatible. 
This condition is similar to those of \cref{thm:4wells}.

Regarding our proof's methods, we combine three ideas. First, the symmetric quasiconvex hull is the set of wells that cannot be separated by symmetric quasiconvex functions. In particular, $-det$ is symmetric quasiconvex. Second, by the characterization of $\mathscr{C}_Q$, we choose suitable translations of quasiconvex functions to identify wells in the convex hull that do not belong to the symmetric quasiconvex hull.  Third, we know explicitely the symmetric lamination convex hull for finite coplanar set explicitly. 

The rest of the paper is organized as follows: for convenience of the reader in \cref{S:section3}, we  briefly quote some results from \cite{camo2020}. This section serves also to clarify and set notation for the rest of the paper. 
In \cref{S:examples} we present explicit examples of well configurations where our results can be applied. The proof of \cref{thm:lamconv} is presented in \cref{S:symmlam}. In Section~\ref{S:young} we give a shorter proof of the \cref{thm:3well} than the one presented in \cite{camo2020}. We devoted the rest of this section to the proof of \cref{thm:4wells} and \cref{thm:iterative}.

\section{The incompatible cone and the three-well problem}\label{S:section3}

In \cite{camo2020}, the authors identify the incompatible cone $\mathcal{C}_U$ of any matrix $U$, see \eqref{eq:cone}. Another equivalent condition for incompatibility between $U$ and $V$ is $\det(V-U)>0$. 
Notice that, if the normal to the plane $\Pi_Q$ is such that $\det Q<0$, then  $\mathcal{C}_U$ (with $U\in\Pi_Q$) intersects $\Pi_Q$ (see Figure~\ref{fig:cones}). Indeed the intersection $\partial \mathcal{C}_U\cap \Pi_Q$ is made of two lines that divide $\Pi_Q$ into its incompatible and compatible components with respect to the vertex of $\mathcal{C}_U$. 

\begin{lemma}\label{planeQV}
Let $Q\in \msym$ such that $\det Q<0$ and $\Pi_Q\subset \re^{2\times 2}_{sym}$ be an affine space with codimension one, such that if $U\in \Pi_Q$, then $\pro{Q}{V-U}=0$ for every $V\in \Pi_Q$. There exist $a$ and $n$ both in $S^1$, not parallel vectors, such that 
\begin{equation}\label{planeQV-1}
\Pi_Q = \{U + \xi\, \ten{a^\perp}{} + \eta\,\ten{n^\perp}{} \,| \, (\xi,\eta)\in \re^2\}.
\end{equation}
Moreover, $W\in \Pi_Q$ and $U$ are compatible if and only if $W\in \mathscr{C}_Q(U).$
\end{lemma}


\begin{remark}\label{det:prop}
In this work we use the convention that  $(\xi,\eta)$ stands for the matrix $\pv{\xi}{\eta}$, where  $a,\,n\in S^1$ and $\xi,\,\eta\in \re$. 
We have that $adj(M) = RMR^T$ for every $M\in \m$, where  $R$ the counter clockwise $\pi/2$-rotation matrix. 
Since $\det M = \pro{adj(M)}{M}$, by direct computation
$\det M = \xi \eta |a\times n|^2$ for every $M\in \Pi_Q$.
\end{remark}

In the symmetric laminar convex hull's characterization for the three-well case when $\det Q<0$, there is an auxiliary matrix $U_0$. This matrix is the vertex of the cone $\mathscr{C}_Q(U_0)$ that defines the boundaries of $L^e(\mathcal{U})$ (see figure~\ref{fig:wedge}). In the next lemma~\cite{camo2020}, we state precisely the conditions for the existence of  $U_0$.  
\begin{lemma}\label{origin}
Let $\{U_1,U_2,U_3\}\subset \msym$ be a linearly independent set of wells such that $\det(U_1-U_2)\leq 0,\ \det(U_1-U_3)\leq 0, \mbox{ and }\det(U_2-U_3)>0.$ Then, there exists $U_0\in C(\mathcal{U})$ such that $\det(U_2-U_0)=\det(U_3-U_0)=0$, and $\det(U_1-U_0)\leq0$.
\end{lemma}

The characterization of $L^e(\mathcal{U})$ for a three-well set in \cite{camo2020} is the following.  
 \begin{proposition}\label{P:lamconv}
Let $ \mathcal{U} = \{U_1,U_2,U_3\}$  be a set of $2\times 2$ linearly independent 
symmetric matrices.
\begin{enumerate}[label=(\alph*)]
\setlength\itemsep{.2em}
\item\label{lhull1}  If $\det(U_1-U_2) > 0$,  $\det(U_1-U_3)> 0$ and $\det(U_2-U_3) \leq 0$, then
$$
L^e(\mathcal{U})=  \{U_1\}  \cup C(\{U_2,U_3\}). 
$$
\item\label{lhull2}  If  $\det(U_1-U_2)\leq 0$, $\det(U_1-U_3)\leq 0$ and $\det(U_2-U_3)> 0$, then
$$
L^e(\mathcal{U})=  C(\{U_0,U_1,U_2\}) \cup C(\{U_0,U_1,U_3\}),
$$
\end{enumerate}
where  $U_0\in C(\mathcal{U})$ is uniquely characterized by  $\det(U_0-U_3)=\det(U_0-U_2)=0$.  
\end{proposition}

\begin{remark}\label{rel-boun-L2} If the elements in the three-well set $\mathcal{U}$ are pairwise compatible, $C(\mathcal{U}) = L^e(\mathcal{U})$, see \cite{Bt,A2016,Z2002}. Indeed, every $U\in C(\mathcal{U})$ is a lamination of degree two. This is readily proved by noticing that the parallel line to the segment $C({U_1,U_2})$ through $U$ is a compatible line and it intersects the segments $C({U_2,U_3})$ and $C({U_3,U_1})$ at two compatible wells, say $V$ and $W$ respectively. Since $L^{e,1}(\mathcal{U})$ is the union of those three segments, we conclude that $U\in L^{e,2}(\mathcal{U})$. Moreover this affirmation holds for every three-well set; this follows directly from the proof of \cref{P:lamconv} in \cite{camo2020}. Therein, the authors prove  that $L^{e,2}(\mathcal{U})=L^{3,2}(\mathcal{U})$.
Therefore, due to the \cref{P:lamconv} and the previous discussion, we conclude that if $\mathcal{U}$ is any set of three wells,  $L^{e,1}(\mathcal{U})\subseteq \rpartial L^{e,2}(\mathcal{U})$. We notice that the equality holds if  either (a) the wells in $\mathcal{U}$ are pairwise incompatible, (b) the wells in $\mathcal{U}$ are pairwise compatible, or (c) there is only one compatible pair of wells in $\mathcal{U}$. In the case where there are two compatibility relations among the elements in $\mathcal{U}$, it follows that $L^{e,1}(\mathcal{U})\subsetneq \rpartial L^{e,2}(\mathcal{U})$.
\end{remark}

In \cite{camo2020}, we also provided the following partial characterization of the quasiconvex hull of $\mathcal{U}$  for a particular family of three-well sets. 
The proof presented in \cite{camo2020} is based on the computation of symmetric polyconvex envelopes for a suitable family of functions. In \cref{S:young} we give an alternative proof. 

\begin{theorem}\label{thm:3well}
Let $\mathcal{U}=\{U_1,U_2,U_3\}$ be a set of $2 \times 2$ linearly independent symmetric wells. If there exist at least two different matrices in $\mathcal{U}$ that are rank-one compatible, then $Q^e(\mathcal{U}) = L^e(\mathcal{U})$. 
\end{theorem}


\section{Explicit examples}\label{S:examples}

\subsubsection*{Four wells example}
Let $\mathcal{U} = \{U_1,\,U_2,\,U_3,\,U_4\}$ where 
\begin{equation}\label{eq:4wells}
U_1=\begin{pmatrix} 1 & 0 \\ 0 & 2 \end{pmatrix}, \quad U_2=\begin{pmatrix} 2 & 0 \\ 0 & -1 \end{pmatrix}, \quad U_3=\begin{pmatrix} -1 & 0 \\ 0 & -2 \end{pmatrix}, \quad \text{and}\quad U_4= \begin{pmatrix} -2 & 0 \\ 0 & 1 \end{pmatrix}.
\end{equation}

\begin{figure}[ht]
     \centering
     \begin{subfigure}[b]{0.4\textwidth}
         \centering
\includegraphics[height=0.2\textheight]{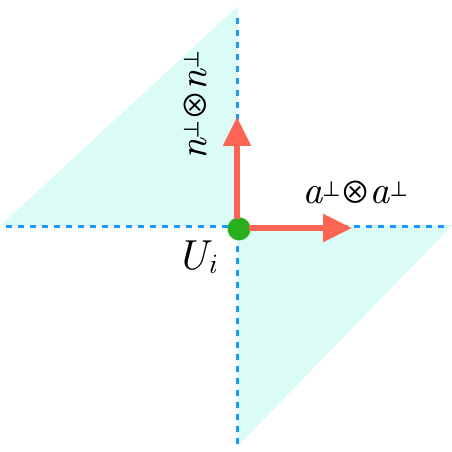}
         \caption{The planar compatible cone $\mathscr{C}_Q(U_i)$}
         \label{fig:example-a}
     \end{subfigure}
     \hspace{1.3cm}
     \begin{subfigure}[b]{0.4\textwidth}
         \centering
         \includegraphics[height=0.2\textheight]{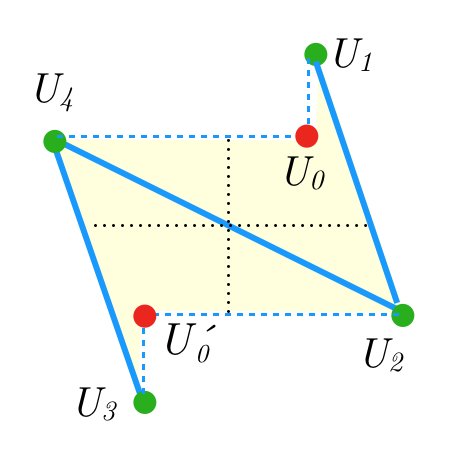}
         \caption{A four-well configuration. The set $L^e(\mathcal\{U_1,\cdots,U_4\})$ is the yellow region.}
         \label{fig:example-b}
     \end{subfigure}
     \caption{The wells $U_1,\cdots,U_4$, are represented with green dots, solid blue lines are compatible segments and dashed lines represent rank-one compatible segments. The auxiliary wells $U_0$ and $U_0'$ are represented with red dots.}
\end{figure}

By a direct computation, we have that $\pro{Q}{U_i}=0$ for every $i=1,2,\cdots, 4$ where
\[Q=\frac{1}{\sqrt{2}}\begin{pmatrix} 0 & 1 \\ 1 & 0 \end{pmatrix} = \sqrt{2}\begin{pmatrix} 1 \\ 0  \end{pmatrix} \odot \begin{pmatrix}  0 \\ 1 \end{pmatrix}.
\]
Hence,  $Q=\nu n\odot a$, with $n=(1,0)^T$ and $a=(0,1)^T$, and the set $\Pi_Q$ contains the origin. The rank-one directions that determine this set are given by, 
\[
n^\perp \odot n^\perp = \begin{pmatrix} 0 & 0 \\ 0 & 1  \end{pmatrix}
, \quad \mbox{and} \quad a^\perp \odot a^\perp = \begin{pmatrix} 1 & 0 \\ 0 & 0  \end{pmatrix}.
\]
Since $\pro{n^\perp \odot n^\perp}{a^\perp \odot a^\perp} = 0$, the rank-one lines spanned by each of these matrices make an angle of $\pi/2$. By \cref{planeQV}, $V\in\Pi_Q$ and $U_i$ are compatible  if and only if $V=U_i + \diag (a,b)$, with $ab \leq 0$, see \cref{fig:example-a}. It follows that $U_2$ is compatible with $U_1$ and $U_4$, also $U_4$ and $U_3$ are also compatible, see \cref{fig:example-b}. Now, let
\[
\mathcal{F} = \left\{ \{U_1,U_2,U_3\},\,  \{U_1,U_3,U_4\},\, \{U_1,U_2,U_4\},\, \{U_2,U_3,U_4\}  \right\}.
\] 
Since inside the first and second sets, only a pair of wells are compatible, from \cref{P:lamconv}, it follows that 
\[
L^e (\{U_1,U_2,U_3\}) = \{U_3\}\cup C(\{U_1,\,U_2\}), \quad \mbox{and} \quad L^e (\{U_1,U_3,U_4\}) = \{U_1\}\cup C(\{U_3,\,U_4\}).
\] 
For the two remaining sets in $\mathcal{F}$, we have two compatibility relations inside each three-well set. 
Hence, by \cref{P:lamconv}, it follows that
\[
\begin{split}
L^e (\{U_1,U_2,U_4\}) = C(\{U_0,\,U_2,\,U_4\})\cup C(\{U_1,\,U_2,\,U_0\}),  \\
L^e (\{U_2,U_3,U_4\}) = C(\{U_0',\,U_2,\,U_4\})\cup C(\{U_4,\,U_3,\,U_0'\}), 
\end{split}
\] 
where $U_0 = \Id\in C(\{U_1,\,U_2,\,U_4\})$ and $U_0' = -\Id\in C(\{U_2,\,U_3,\,U_4\})$ are the solutions of  
\[ \det(U_1-U_0) = \det(U_4-U_0) = 0, \quad \mbox{and} \quad \det(U_3-U_0')=\det(U_2-U_0')=0.
\] 
Now, we notice that there are only one compatible pair in $\{U_1,U_2,U_3\}$ and $\{U_1,U_3,U_4\}$ within each set, but  $\{U_1,U_2,U_4\}$ and  $\{U_2,U_3,U_4\}$, have two compatible pairs within each set. 
Therefore, $L^{e,1} (\{U_1,U_2,U_3\}) \cup L^{e,1} (\{U_1,U_3,U_4\}) \subset L^{e,1} (\{U_1,U_2,U_4\}) \cup L^{e,1} (\{U_2,U_3,U_4\})$ and in turns, 
$L^e (\{U_1,U_2,U_3\}) \cup L^e (\{U_1,U_3,U_4\}) \subset L^e (\{U_1,U_2,U_4\}) \cup L^e (\{U_2,U_3,U_4\})$.
Hence, by \cref{thm:lamconv}, the symmetric lamination convex hull of the set $\mathcal{U}$, see \cref{fig:example-b}, is 
\begin{equation}\label{eq:4well-lam}
L^e(\mathcal{U}) = C(\{U_0',\,U_2,\,U_4\})\cup C(\{U_4,\,U_3,\,U_0'\}) \cup C(\{U_0',\,U_2,\,U_4\})\cup C(\{U_4,\,U_3,\,U_0'\}).
\end{equation}

\subsubsection*{A degenerated four-well example}

We consider a {\em degenerate} case of the previous four-well problem.
Let $\mathcal{U} = \{U_1,\,U_2',\,U_3,\,U_4'\}$ such that $U_1$ and $U_3$ are as in \cref{eq:4wells} and 
\[
U_2'=\begin{pmatrix} 1 & 0 \\ 0 & -1 \end{pmatrix}, \quad \mbox{and} \quad U_4'= \begin{pmatrix} -1 & 0 \\ 0 & 1 \end{pmatrix}.
\] 
\begin{figure}[ht]
\begin{centering}
\includegraphics[scale=0.5]{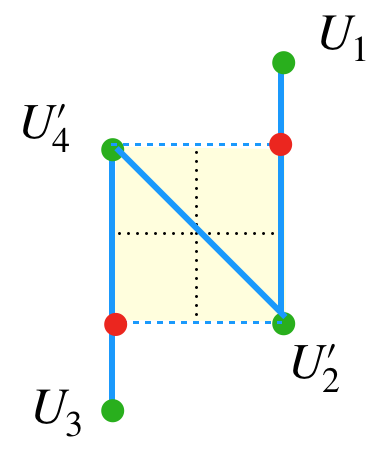}
\caption{A four-well configuration where two pairs of wells are rank-one compatible. Here $Q^e(\mathcal{U})=L^e(\mathcal{U})$ and both sets are strictly contained in $C(\mathcal{U})$
}
\label{fig:2flags}
\end{centering}
\end{figure}
In this case $Q$, $a,\, n,\, U_0$ and $U_0'$ are as in the previous example. 
Thus, $L^e(\mathcal{U})$ is as in \cref{eq:4well-lam} with $U_2$ and $U_4$ replaced by $U_2'$ and $U_4'$, respectively. 
This four-well set satisfies the conditions of  \cref{it:le1} in  \cref{thm:4wells} and  $Q^e(\mathcal{U})=L^e(\mathcal{U})$. Notice that this configuration is the union of two flag configurations, see \cref{fig:2flags}.

\subsubsection*{Five wells example}
Now we consider a five-well problem. Let $\mathcal{U}=\{U_1, \cdots,\,U_5\}$ where 
\begin{equation}\label{eq:5wells}
\begin{array}{c}
U_1=\begin{pmatrix} 0 & 0 \\ 0 & 2 \end{pmatrix}, \quad
U_2=\begin{pmatrix} 1 & 1 \\ 1 & 5 \end{pmatrix}, \quad
U_3=\begin{pmatrix} 2 & 2 \\ 2 & 6 \end{pmatrix},\\ \\
U_4=\begin{pmatrix} 3 & 3 \\ 3 & 3 \end{pmatrix}, \quad
\text{and}\quad
U_5= \begin{pmatrix} 2 & 2 \\ 2 & 0 \end{pmatrix}
\end{array}
\end{equation}
By direct computation, we have that $\pro{Q}{U_i}=0$ for every $i=1,2,3,4$, where 
\[
Q=\frac{1}{\sqrt{6}}\begin{pmatrix} 2&-1\\-1&0 \end{pmatrix}.
\]
Hence, $\mathcal{U}$ is coplanar and the plane $\Pi_Q$ is well defined. Due to Lemma~1 in \cite{camo2020} and $\det Q< 0$, we have that $Q=\nu a\odot n$, 
where $\nu$ is  a real number and $a,n\in S^1$. 
Notice that these vectors are linear combinations of the eigenvectors of $Q$, thus 
\begin{equation}
a=\frac{1}{\sqrt{2}}\begin{pmatrix}1\\-1
\end{pmatrix},\ \ n=\begin{pmatrix}1\\ 0
\end{pmatrix} \quad \mbox{and} \quad a^\perp\otimes a^\perp = \frac{1}{2}\begin{pmatrix}1 & 1\\ 1 & 1 \end{pmatrix}, \ \  n^\perp\otimes n^\perp = \begin{pmatrix} 0 & 0 \\ 0 & 1 \end{pmatrix}.
\end{equation}

\begin{figure}[ht]
     \centering
     \begin{subfigure}[b]{0.4\textwidth}
         \centering
\includegraphics[height=0.2\textheight]{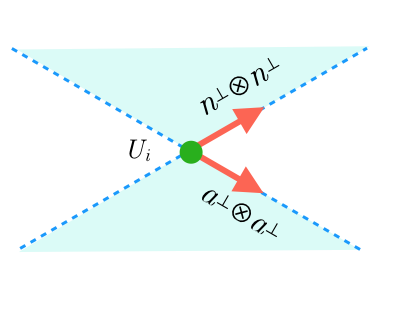}
         \caption{The planar compatible cone $\mathscr{C}_Q(U_i)$.}
         \label{fig:5wells-a}
     \end{subfigure}
     \hspace{1.3cm}
     \begin{subfigure}[b]{0.4\textwidth}
         \centering
         \includegraphics[height=0.2\textheight]{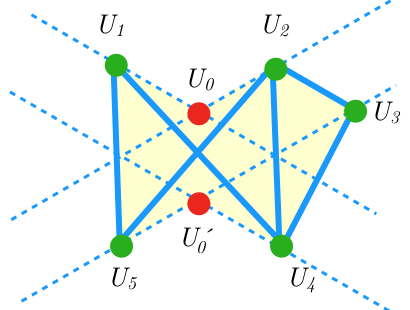}
         \caption{A five-well configuration. The set $L^e(\mathcal\{U_1,\cdots,U_4\})$ is the yellow region.}
         \label{fig:5wells-b}
     \end{subfigure}
     \caption{The five-well example. The wells $U_1,\cdots,U_5$, are represented with green dots, solid blue lines are compatible segments and dashed lines represent rank-one compatible segments. The auxiliary wells $U_0$ and $U_0'$ are represented with red dots.}

\end{figure}

In this case $\pro{a^\perp\otimes a^\perp}{n^\perp\otimes n^\perp}=1/2$, and the angle between the rank-one 
lines generated by those matrices is $\pi/6$, see \cref{fig:5wells-a}. 
Arguing as in the previous example, we conclude that
\[
L^e(\mathcal{U}) =L^e(\{U_1,U_2,U_5\})\cup L^e(\{U_2,U_3,U_4\}) \cup L^e(\{U_2,U_4,U_5\}).
\]
Moreover, by the compatibility relations among the wells, see \cref{fig:5wells-b}, we get
\[
\begin{array}{c}
L^e(\{U_1,U_2,U_5\}) = C(\{U_0,U_1,U_5\})\cup C(\{U_0,U_5,U_2\}),\\
L^e(\{U_2,U_4,U_5\}) = C(\{U_2,U_4,U_0'\})\cup C(\{U_2,U_0',U_5\}),\\
L^e(\{U_2,U_3,U_4\}) = C(\{U_2,U_3,U_4\}),
\end{array}
\]
where $U_0\in C(\{U_1,U_2,U_5\})$ and $U_0'\in C(\{U_2,U_4,U_5\})$ are the solutions of (see \cref{fig:5wells-b}),  
\[
\det{(U_1-U_0)}=\det{(U_2-U_0)}=0, \quad \mbox{and}\quad \det{(U_4-U_0')}=\det{(U_5-U_0')}=0.
\]

\section{The symmetric lamination convex hull}
\label{S:symmlam}

This section is devoted to prove \cref{thm:lamconv}. To this end, we show the reverse inclusion in
\cref{laminarcontation}, namely
\begin{equation*}
L^{e,i}(\mathcal{U}) \subset  \bigcup_{\mathcal{V}\in \mathcal{F}}L^{e,i}(\mathcal{V}).
\end{equation*}

\subsection{Preliminary lemmas}

Along this subsection we will assume that each set of coplanar wells is contained in the affine space $\Pi_Q$ where its normal matrix $Q$ has negative determinant.

\begin{lemma}\label{supcone}
Let $U,V$ and $W\in \re^{2\times 2}_{sym}$ such that $V\in\mathscr{C}^+_Q(U)$ and $W\in\mathscr{C}^-_Q(U)$, then $V\in\mathscr{C}^+_Q(W)$. 
\end{lemma}
\begin{proof}
Let $a$ and $n$ be the unitary vectors associated with $Q$ as in the definition of $\mathscr{C}_Q(U)$. The result is  straight forward, since
\[
V = U + (\xi, \eta), \text{ with } \xi\geq 0\geq \eta\quad\text{ and }\quad 
W = U + (\gamma, \delta), \text{ with }  \delta\geq 0\geq \gamma,
\]
imply that $V=W + (\xi-\gamma, \eta-\delta)$, where $ \xi-\gamma\geq 0\geq  \eta-\delta$ and we have 
used the notation introduced in \cref{det:prop}.
\end{proof}

We say that the matrices of a coplanar set $\mathcal{U}=\{U_1, \cdots,\,U_n\}$ are labeled in increasing order if at any edge of $C(U)$ the corresponding vertices are $U_i$ and 
$U_{i+1}$ for $i\in\{1,\cdots,n\}$ in cyclic order.


\begin{lemma} \label{lem:Diag}Let $\mathcal{L}:=\{U_1,U_2,U_3,U_4\}\subset\re^{2\times 2}_{sym}$ be a coplanar set such that its convex hull $C(\mathcal{L})$ is a quadrilateral labeled in increasing order.
If there are three values for  $i\in\{1,2,3,4\}$ such that the pairs $\{U_i,U_{i+1}\}$  are compatible, then either $\{U_1,U_3\}$ or $\{U_2,U_4\}$ are compatible.   
\end{lemma}
\begin{proof}
Without loss of generality,  we  assume $U_{1}=0$. By Lemma~1 in \cite{camo2020} and \cref{det:prop}, there exist $b,c,d,f,g,h \in \re$ such that 
\[
U_{1}=(0,0), \quad U_{2}=(b,c),\quad U_{3}=(d,f),\quad \mbox{and} \quad U_{4}=(g,h). 
\]
 Without loss of generality, we assume that for each $i\in\{1,2,3\}$ the wells $U_i,U_{i+1}$ are compatible but $U_{1},\,U_{4}$ are incompatible wells. Hence,  
we get (up to a common factor $|a\times n|^2$) the following set of equations
\begin{equation} \label{eq:rels1}
0\geq \det (U_{2}-U_{1}) = bc,
\end{equation}
\begin{equation} \label{eq:rels2}
0\geq \det (U_{3}-U_{2})=(d-b)(f-c),
\end{equation}
\begin{equation} \label{eq:rels3}
0\geq \det (U_{4}-U_{3})=(g-d)(h-f),
\end{equation}
\begin{equation} \label{eq:rels4}
gh=\det (U_{4}-U_{1})>0.
\end{equation}
From \cref{eq:rels1} we have two options. First, we assume that $c\leq 0 \leq b$.
Due to \cref{eq:rels2}, we have two more cases, either 
\[
\text{(I) } (f-c)\leq 0\leq (d-b) \quad\text{ or }\quad 
\text{(II) } (d-b)\leq 0\leq (f-c).
\]
If (I) holds, it follows that $U_{3}\in \mathscr{C}^+_{Q}(U_{2})$ and $U_{1}\in \mathscr{C}^-_{Q}(U_{2})$. 
Then, by \cref{supcone}, the wells $U_{1}$ and $U_{3}$ are compatible.
Now, if (II) holds, then $U_{2} = U_{3} +(b-d, c-f)$, and $U_{2}\in \mathscr{C}^+_{Q}(U_{3})$.
We notice that $U_{3}$ satisfies either $0< fb-cd$ and the set $\mathcal{L}$ is clockwise oriented, or $ fb-cd < 0$ and the set $\mathcal{L}$ is counterclockwise oriented.
To preserve the vertex orientation and relations \eqref{eq:rels3} and \eqref{eq:rels4}, it follows that $U_{4}\in\mathscr{C}^-_{Q}(U_{3})$ in either case. 
Hence, by \cref{supcone}, the wells $U_{4}$ and $U_{2}$ are compatible.   
 
Second, we assume $b\leq 0 \leq c$. The arguments in this case follows the same lines as the previous case 
and we skip the proof.
\end{proof}
\begin{figure}[ht]
\begin{centering}
\includegraphics[scale=0.5]{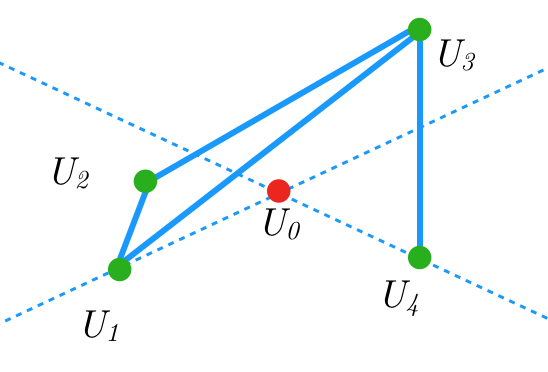}
\caption{A four-well configuration where one the boundary segments is not a compatible one.}
\label{4-wells-order}
\end{centering}
\end{figure}
\begin{remark}\label{rmk:det}
It is well known, see \cite{BFJK,BDa}, that if $f:\re^{m\times m}\rightarrow \re$ is quadratic and $f(a\otimes n)\geq 0$ for every $a,n\in\re^m$ then $f$ is rank-one convex. Moreover, if both conditions are satisfied by the function $f$, then it is quasiconvex. This result and \cref{det:prop} easily imply that $-\det(\cdot):\re^{2\times 2}_{sym}\rightarrow \re$ is a symmetric quasiconvex function.
\end{remark}

\begin{lemma}\label{lem:linedot}
Let $U,V,W\in \re^{2\times 2}_{sym}$ such that $V,W$ are compatible and  $U$ is incompatible with both of them. Then $U$ is incompatible with any point in $C(\{U,V,W\})\setminus \{U\}$.  
 \end{lemma}
 
\begin{proof}
Let $M\in C(\{U,V,W\})\setminus \{U\}$. 
Hence, $M=\lambda_1U + \lambda_2V + \lambda_3W$ for some $\lambda_1\in[0,1)$ and $\lambda_2,\lambda_3\in [0,1]$,  such that $\lambda_1+\lambda_2+\lambda_3 =1$.
Since $\lambda_1\neq 1$, we have that $\lambda_2 + \lambda_3>0$ and 
$$
M-U =(\lambda_2+\lambda_3)\left(s(V-U) + (1-s)(W-U)\right),
$$
where $s=\lambda_2/(\lambda_2+\lambda_3)$, and $(1-s)= \lambda_3/(\lambda_2+\lambda_3)$.
Hence, $\det(M-U) = (\lambda_2+\lambda_3)^2 \det(s(V-U) + (1-s)(W-U))$. 
Now, since $-\det(e(\cdot))$ is a rank-one convex function, and $\{(V-U),(W-U)\}$ is a  compatible set, we have that
\[
\det(s(V-U) +(1-s)(W-U)) \geq s\det(V-U)+(1-s)\det(W-U)>0.
\]
Therefore, $U$ and $M$ are incompatible and the proof is complete.
\end{proof}

 \begin{lemma}\label{lem:4wells-2}
Let $\mathcal{L}:=\{U_{v,1},U_{v,2},U_{w,1} U_{w,2}\}\subset \re^{2\times 2}_{sym}$ be coplanar four-well set with $\det(U_{v,1}-U_{v,2})\leq 0$ and $\det(U_{w,1}-U_{w,2})\leq 0$. Also let $\mathcal{F}$ as in \cref{laminarcontation} and  assume there exist $u\in C(\mathcal{L})$, $v\in L^{e,1}(\{U_{v,1},U_{v,2}\})$, and $w\in L^{e,1}(\{U_{w,1},U_{w,2}\})$ such that  $u\in L^{e,1}(\{v,w\})$. Then, $u\in L^{e,2}(\mathcal{V})$ for some $\mathcal{V}\in \mathcal{F}$.
 \end{lemma}
\begin{proof}
We consider three different cases depending on the number of extreme points (or vertices) of $C({\mathcal{L}})$.

First, if the set of extreme points of $C({\mathcal{L}})$ consists of two points, $\mathcal{L}$ is contained in a compatible line, and the affirmation follows trivially. 
Second, we assume that $C(\{U_{v,1},U_{v,2},U_{w,1},U_{w,2}\})$ has exactly three extremal points. 
Without loss of generality, let 
\begin{equation}\label{eq:triangle}
U_{w,2}\in C(\{U_{v,1},U_{v,2},U_{w,1}\}).
\end{equation}
By \cref{lem:linedot}, either $U_{v,1}$ or $U_{v,2}$  are compatible with $U_{w,1}$, otherwise $U_{w,2}$ would be incompatible with $U_{w,1}$. Hence, there exist at least two compatibility relations between the extreme points  $U_{v,1},U_{v,2},$ and $U_{w,1}$. By \cref{P:lamconv}, since $U_{w,2}$ and $U_{w,1}$ are compatible, the whole segment $\{tU_{w,1} + (1-t)U_{w,2}|t\in[0,1]\}$ is contained in $L^{e,2}(\{U_{v,1},U_{v,2},U_{w,1}\})$. Therefore $v,w\in L^{e,2}(\{U_{v,1},U_{v,2},U_{w,1}\})$ and $u\in L^{e,3}(\{U_{v,1},U_{v,2},U_{w,1}\})$. Due to \cref{P:lamconv} and \cref{rel-boun-L2}, $L^{e,2}(\{U_{v,1},U_{v,2},U_{w,1}\}) = L^{e}(\{U_{v,1},U_{v,2},U_{w,1}\})$. Hence, $u\in L^{e,2}(\{U_{v,1},U_{v,2},U_{w,1}\})$ and the lemma follows. 
\begin{figure}[ht]
\begin{centering} 
\includegraphics[scale=0.55]{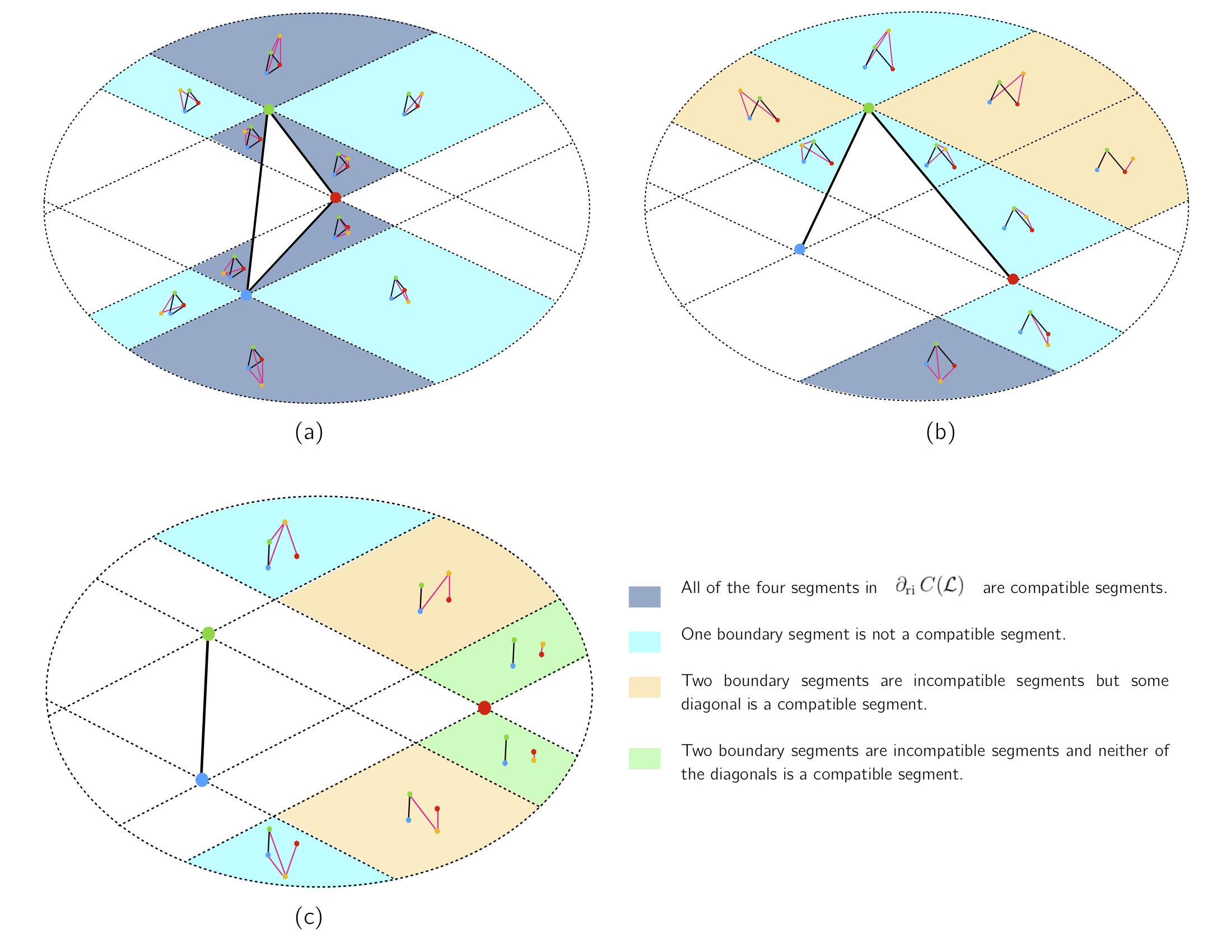}
\caption{Figures (a), (b) and (c) give different sections of the plane $\Pi_Q$ where $U_{w,2}$ can be placed. The positions for $U_{v,1},\, U_{v,2}$ and $U_{w,1}$ are fix in each figure but the number of compatibility relations among these wells is different. Dotted lines represent boundaries of the compatibility cone centered at different points. In figure (a), we assume that $U_{v,1},\, U_{v,2}$ and $U_{w,1}$ are pairwise compatible. In Figure (b) $U_{v,1}$ and $U_{w,1}$ are incompatible. In Figure (c) $U_{w,1}$ is incompatible with $U_{v,1}$ and $U_{v,2}$.
}
\label{fig:config}
\end{centering}
\end{figure}

Third, we assume that all wells in  $\mathcal{L}$ are the extremal points of $C(\mathcal{L})$. 
Since $\mathcal{L}$ is coplanar, the set $C(\mathcal{L})$ is a quadrilateral, and there exist six segments joining the wells in $\mathcal{L}$. 
The boundary $\rpartial C(\mathcal{L})$ is given by 
four segments on the boundary of the quadrilateral.
The two remaining segments connect the opposite wells of $\mathcal{L}$ diagonally. 
Now, we split our analysis into three cases that depend on how many compatible segments are on $\rpartial C(\mathcal{L})$, see \cref{fig:config}.

\begin{enumerate}[label=\alph*)]
\item \label{step-1}
If all four segments in  $\rpartial C(\mathcal{L})$ are compatible, then either $U_{w,1}$ and $U_{v,2}$ or $U_{w,2}$ and $U_{v,1}$  are compatible. This follows from the assignment $(U_{v,1},U_{v,2},U_{w,2},U_{w,1})\mapsto (U_1,U_2,U_3,U_4)$ and \cref{lem:Diag}. 
Hence, a diagonal in $C(\mathcal{L})$ is a compatible segment, and every point in $C(\mathcal{L})$ belongs to a triangle with pairwise compatible vertices and $u\in L^{e,2}(\mathcal{V})$ for some $\mathcal{V}\in \mathcal{F}$.

\item \label{step-2} We assume that only one out of the four segments in $\rpartial C(\mathcal{L})$ is not a compatible segment. 
By assignment $(U_{v,1},U_{v,2},U_{w,2},U_{w,1})\mapsto (U_1,U_2,U_3,U_4)$ and \cref{lem:Diag}, at least one of the diagonal segments is compatible.
Thus, $C(\mathcal{L}) = T_1\cup T_2$, where $T_1$ and $T_2$ are two triangles such that $\rpartial(T_1)$ and $\rpartial(T_2)$ consist of three and only two compatible segments, respectively. 
If $u\in T_1$ then $u\in L^{e,2}(\mathcal{V})$ for some $\mathcal{V}\in \mathcal{F}$ and the lemma follows. 
Hence, let $u\in T_2$ and, without loss of generality, assume that the set of extreme points of $T_2$ is ${\cal L}_T=\{U_{v,1},\, U_{w,1},\, U_{w,2}\}$ with the incompatible line segment $\{tU_{v,1} + (1-t)U_{w,1}\,|\, t\in [0,1]\}$.
Notice that $u$ is a convex combination of two compatible wells $v\in L^{e,1}(\{U_{v,1},U_{v,2}\})$ and $w\in L^{e,1}(\{U_{w,1},U_{w,2}\})$, hence  there exists $v'\in L^{e,1}(\{U_{v,1},U_{w,2}\})$ such that $u$ also is a convex combination of $v'$ and $w$, due to $\mathcal{L}_T\subset \mathcal{L}$. Moreover $v'$ and $w$ are compatible since they belong to the compatible segment $\{tv + (1-t)w\,|\, t\in [0,1]\}$, so we conclude that $u\in L^{e,2}(\mathcal{L}_T)$.


\item \label{step-3}
We assume there are only two compatible segments on the boundary of $C({\cal L})$. 
We claim that there exists a compatible diagonal segment.  
If the claim holds, we let $\{tU_{v,1} + (1-t)U_{w,2}\,|\, t\in [0,1]\}$ be the compatible diagonal segment without loss of generality. 
Hence, $C(\mathcal{L}) = T_1\cup T_2$ where $T_1 = C(\{U_{v,1},\, U_{w,2},\, U_{v,2}\})$ and $T_2 = C(\{U_{v,1},\, U_{w,1},\, U_{w,2}\}) $ are triangles such that $\rpartial T_1$ and $\rpartial T_2$ have two compatible segments, see \cref{fig:4well1}. 
We have that either $u\in T_1$ or $u\in T_2$. 
In both cases, we argue as in the second part of \ref{step-2} to prove that $u\in L^{e,2}({\cal V})$ 
where either ${\cal V} =\{U_{v,1},\, U_{w,1},\, U_{w,2}\}$ or  ${\cal V} = \{U_{v,1},\, U_{w,2},\, U_{v,2}\}$ and the lemma follows. 
\end{enumerate}

\begin{figure}[ht]
     \centering
     \begin{subfigure}[b]{0.3\textwidth}
         \centering
         \includegraphics[width=\textwidth]{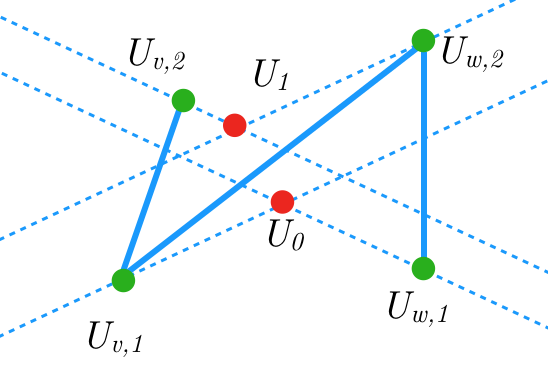}
         \caption{$\ $}
         \label{fig:4well1}
     \end{subfigure}
     \hspace{2cm}
     \begin{subfigure}[b]{0.3\textwidth}
         \centering
         \includegraphics[width=\textwidth]{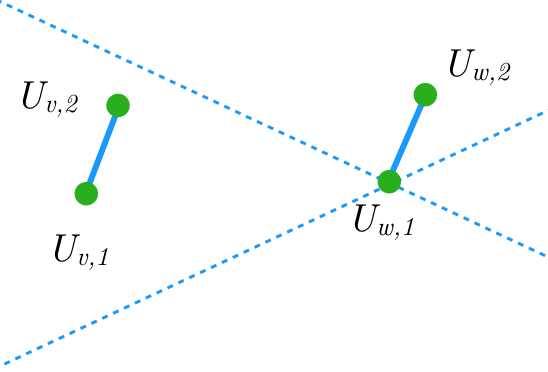}
         \caption{$\ $}
         \label{fig:4well2}
     \end{subfigure}
     \caption{Figure (a) and (b) shows a four-well configuration  where two of the boundary segments are incompatible segments. Figure (b) is not considered in the statement of \cref{lem:4wells-2} since there is not a well $C(\mathcal{L})$ that meets the conditions of \cref{lem:4wells-2}.}
\end{figure}

Finally, we prove the claim in part \ref{step-3}. By contradiction, we assume that there are no compatible diagonals segments, namely $\{U_{w,1},U_{v,2}\}$ and $\{U_{w,2},U_{v,1}\}$ are incompatible. By \cref{lem:linedot}, it follows that $U_{v,1}$ and $U_{v,2}$ are incompatible with every element in $\{sU_{w,1} +(1-s)U_{w,2}\,|\, s\in [0,1]\}$. Hence, on the one hand, for every $s\in[0,1]$, we have that
$$
\det(sU_{w,1} +(1-s)U_{w,2}-U_{v,1})>0\text{ and }
\det(sU_{w,1} +(1-s)U_{w,2}-U_{v,2})>0.
$$
On the other hand, if $t\in [0,1]$, the rank-one convexity of $-\det(e(\cdot))$ and the compatibility of $U_{v,1}$ and $U_{v,2}$ implies that 
\[
\begin{split}
\det(w-v) &= \det(sU_{w,1} +(1-s)U_{w,2}-tU_{v,1} -(1-t)U_{v,2})\\ 
& =\det[t(sU_{w,1} +(1-s)U_{w,2}-U_{v,1})+ (1-t)(sU_{w,1} +(1-s)U_{w,2}-U_{v,2})]\\
& \geq t\det(sU_{w,1} +(1-s)U_{w,2}-U_{v,1}) + (1-t)\det(sU_{w,1} +(1-s)U_{w,2}-U_{v,2})\\ 
& > 0.
\end{split}
\]
Since by assumption there exists $u\in L^{e,1}(\{v,w\})$ and $v$ and $w$ are compatible,  $\det(v-w)\le 0$, a contradiction, and we finish the proof.

\end{proof}

 \begin{lemma}\label{lem:4wells}
Let $\mathcal{U}\subset \re^{2\times 2}_{sym}$ be a set of n coplanar wells and $\mathcal{F}$ as in \cref{laminarcontation}.  Also let $\mathcal{V}_1,\,\mathcal{V}_2$ be two different sets in $\mathcal{F}$.  If there exists $u\in\re^{2\times 2}_{sym}$ such that $u$ is a symmetric lamination of degree one of two compatible wells, $v\in L^{e,2}(\mathcal{V}_1)$ and $w\in L^{e,2}(\mathcal{V}_2)$. 
Then, $u$ is a symmetric lamination of degree one of two compatible wells, $v'\in L^{e,1}(\mathcal{V}_1)$ and  $w'\in L^{e,1}(\mathcal{V}_2)$.
\end{lemma}
\begin{proof}
We notice that if $u\in L^{e,2}(\mathcal{V}_1)\cup L^{e,2}(\mathcal{V}_2)$ there is nothing  to prove, hence we focus on the case $u\notin L^{e,2}(\mathcal{V}_1)\cup L^{e,2}(\mathcal{V}_2)$. 

We claim that there exists $w'\in L^{e,1}(\mathcal{V}_2)$ such that $u\in L^{e,1}(\{v,w'\})$. 
Notice that by the assignment $v\mapsto w$ and $w'\mapsto v$, the claim implies that there exists $v'\in L^{e,1}(\mathcal{V}_1)$ such that $u\in L^{e,1}(\{v',w'\})$ and 
the result follows.

Now we prove the claim. If $w\in L^{e,1}(\mathcal{V}_2)$, by letting $w'=w$ the affirmation follows straight forward. 
Hence, we assume that $w\in L^{e,2}(\mathcal{V}_2)\setminus L^{e,1}(\mathcal{V}_2)$. 
Let $\ell_0,\ell_1$ and $\ell_2$ denote line segments that make up the boundary $\rpartial C(\mathcal{V}_2)$ and  $\mathcal{L}$ be the compatible line through $v$ and $w$. Since $w$ and $v$ belong to the relative interior and exterior of $C(\mathcal{V}_2)$, respectively, $\mathcal{L}$ intersects two of the line segments $\ell_0,\ell_1$ and $\ell_2$. 
Now, by \cref{rel-boun-L2} and \cref{P:lamconv} it follows that
\[
L^{e,1}(\mathcal{V}_2) = \begin{cases} \rpartial C(\mathcal{V}_2) & \mbox{if the wells in } \mathcal{V}_2 \mbox{ are pairwise compatible}\\
\ell_i\cup \ell_{i+1}, \mbox{for some } i=1,2,3 & \mbox{if there is a incompatible pair in } \mathcal{V}_2.
\end{cases}
\]
Therefore, $\mathcal{L}\cap L^{e,1}(\mathcal{V}_2)\neq \emptyset$, and the claim follows by choosing 
$$w'=\mbox{arg}\max\left\{ \mbox{dist} (v,z)\,|\, z\in L^{e,1}(\mathcal{V}_2) \cap \mathcal{L}\right\}.$$ 
 \end{proof}
 
 
 \subsection{The symmetric lamination convex hull for the coplanar \texorpdfstring{$n$}{}-well problem}
 
\begin{proof}[Theorem~\ref{thm:lamconv}]
By definition of $L^{e,1}(\mathcal{U})$ and the equality $L^{e,0}(\mathcal{U})=\mathcal{U}$, it follows that 
\begin{equation}\label{eq:thm1-1}
    L^{e,1}(\mathcal{U}) = \cup_{\mathcal{V}\in \mathcal{F}}L^{e,1}(\mathcal{V}).
\end{equation}
Now, we claim that
\[
\bigcup_{\mathcal{V}\in \mathcal{F}}L^{e,2}(\mathcal{V}) = L^{e,2}(\mathcal{U}). 
\]
By contradiction, we assume that there exists 
\begin{equation}\label{eq:hypothesis1}
u\in L^{e,2}(\mathcal{U})\setminus  \cup_{\mathcal{V}\in \mathcal{F}}L^{e,2}(\mathcal{V}).
\end{equation}
So $u$ is a convex combination of two elements in the set of laminations
of degree one $v,w \in L^{e,1}(\mathcal{U})$. Due to \eqref{eq:thm1-1} $v\in L^{e,1}(\mathcal{V}_1)$ and $w\in L^{e,1}(\mathcal{V}_2)$ for some $\mathcal{V}_1,\mathcal{V}_2\in 
\mathcal{F}$.
We assume that  $\mathcal{V}_1\neq \mathcal{V}_2$, otherwise $u$ would belong to 
$L^{e,2}(\mathcal{V}_1)$, a contradiction. 
Moreover, if $v\in \mathcal{V}_1$ or $w\in \mathcal{V}_2$, 
then there exists $\mathcal{V}\in\mathcal{F} $ such that $u \in  
L^{e,2}(\mathcal{V})$, again a contradiction. 
Hence, there are four wells 
$U_{v,1},U_{v,2}\in \mathcal{V}_1$ and  $U_{w,1},U_{w,2}\in \mathcal{V}_2$ such that 
$v$ and $w$ belong to the compatible segments $L^{e,1}(\{U_{v,1},U_{v,2}\})$ and 
$L^{e,1}(\{U_{w,1},U_{w,2}\})$, respectively, and $u\in C(\{U_{v,1},U_{v,2},U_{w,1},U_{w,2}\})$. 
By \cref{lem:4wells-2}, $u\in L^{e,2}(\mathcal{V})$ for some $\mathcal{V}\in \mathcal{F}$, a contradiction, and the claim follows.
Thus, by \cref{rel-boun-L2} 
$L^{e,2}(\mathcal{V})=L^e(\mathcal{V})$ for every three-well set $\nu$ and 
\begin{equation}\label{lamination2}
L^{e,2}(\mathcal{U})=\bigcup_{\mathcal{V}\in \mathcal{F}}L^{e,2}(\mathcal{V})=
\bigcup_{\mathcal{V}
\in \mathcal{F}}L^{e}(\mathcal{V}). 
\end{equation}

Next, prove that 
\[
L^{e,3}(\mathcal{U})=L^{e,2}(\mathcal{U}).
\]
We proceed again by contradiction. Assume that $u \in L^{e,3}(\mathcal{U})\setminus L^{e,2}(\mathcal{U})\neq \emptyset$, so it is a convex combination of two compatible wells, say   $v,w\in L^{e,2}(\mathcal{U})$,  and by the first equality in \cref{lamination2},  $v\in L^{e,2}(\mathcal{V}_1)$ and $w\in L^{e,2}(\mathcal{V}_2)$ for some $\mathcal{V}_1,\, \mathcal{V}_2 \in \mathcal{F}$. Due to the last equality in \cref{lamination2} we also have that 
$\mathcal{V}_1\neq \mathcal{V}_2$. Thus, by \cref{lem:4wells}, $u$ is a convex combination of two compatible wells 
$v'\in L^{e,1}(\mathcal{V}_1)$ and $w'\in L^{e,1}(\mathcal{V}_1)$, where $\mathcal{V}_1\neq \mathcal{V}_2$.
So, there exist two pair of compatible wells $U_{v,1},U_{v,2}\in \mathcal{V}_1$ and $U_{w,1},U_{w,2}\in \mathcal{V}_1$ such that $v'$ and $w'$ are convex combinations of the former and latter pairs, respectively, and by \cref{lem:4wells-2}, $u\in L^{e,2}(\mathcal{V})$ for some $\mathcal{V}\in \mathcal{F}$, a contradiction. Therefore $L^{e,3}(\mathcal{U})= L^{e,2}(\mathcal{U})$ and the proof is completed.
\end{proof}

\section{The quasiconvex hull for the coplanar \texorpdfstring{$n$}{}-well problem. }\label{S:young}

The symmetric quasiconvex hull $Q^e(\mathcal{U})$ is the set of wells that cannot be separated from the set $\mathcal{U}$ by symmetric quasiconvex functions, see \cref{semihulls}. In the next lemma, we provide sufficient conditions on the elements in $C(\mathcal{U})$ that do not belong to $Q^e(\mathcal{U})$.

\begin{lemma}\label{quasidet}
Let $\mathcal{U}\subset \re^{2\times 2}_{sym}$ be a finite coplanar set of wells contained in the affine space $\Pi_Q$. Assume there exists $U_0\in \Pi_Q$  such that $\det(V-U_0)\geq 0$ for every $V\in \mathcal{U}$. If $U\in C(\mathcal{U})$ and $\det(U-U_0)<0$, then $U\notin Q^e(\mathcal{U})$. 
\end{lemma}

\begin{proof} 
We proceed by contradiction and assume that $U\in Q^e(\mathcal{U})$.
From \cref{rmk:det}, the function $-\det :\re^{2\times 2}_{sym}\rightarrow \re$ 
and its translation $f(\cdot)=-\det((\cdot) -U_0)$ are both symmetric quasiconvex. 
Due to \cref{semihulls}, $U$ satisfies $-\det (U-U_0) \leq \sup\{ -\det (V-U_0)\,\,|\, V\in \mathcal{U}\}$, but this is a contradiction to the statement, so $U\notin Q^e(\mathcal{U})$.
\end{proof}

\subsection{The three-well case}

Theorem~\ref{thm:3well} is already proven in \cite{camo2020}.
Here, we present an alternative shorter proof.

\begin{proof}[Theorem~\ref{thm:3well}]
We divide the proof into three cases corresponding to the number of pairwise compatible pairs formed by the wells in $\mathcal{U}$.

\begin{enumerate}[label=(\alph*)]
\item In the case of three compatible pairs, Bhattacharya et al. proved in \cite{Bt} that  
$Q^e(\mathcal{U})=L^e(\mathcal{U})=C(\mathcal{U})$. The affirmation follows as a particular case.

\item 
We assume that only two distinct pairs in $\mathcal{U}$ are compatible, and one of them is rank one compatible. Without loss of generality, let 
\begin{equation*}
\det(U_1-U_3)\le 0, \det(U_1-U_2)=0 \text{ and } \det(U_2-U_3)>0. 
\end{equation*}
Now, we claim  there exists $U_0\in C(\{U_1,U_2\})$ such that 
\begin{equation}\label{2conect}
\det(U_0-U_1)=\det (U_0-U_2)=\det (U_0-U_3)=0.
\end{equation}
Indeed, if $\det(U_1-U_3)=0$, we let $U_0 = U_1$ and \eqref{2conect} follows. If $\det(U_1-U_3)<0$, by the continuity of  $t\mapsto \det(tU_1+(1-t)U_2-U_3)$ there exists a  $t_0\in (0,1)$ such that $U_0 = t_0U_1+(1-t_0)U_2$
satisfies the necessary condition and the claim follows.

By the assumptions on $\mathcal{U}$, the matrices $U_2-U_0$ and $U_3-U_0$ are linearly independent. 
Without losing generality, we let $U_2 = U_0+(\alpha,0)$ and $U_3 = U_0+(0,\beta)$, where $\alpha\beta<0$ due to the incompatibility of $U_2$ and $U_3$. Now, let $U\notin C(\mathcal{U})\setminus L^e(\mathcal{U})$. 
As a consequence of \cref{P:lamconv}, $U\in \relint C(\{U_0,U_2,U_3\})\cup C(\{U_2,U_3\})$, and $U = \lambda_0U_0+\lambda_2U_2+\lambda_3U_3$ for some  $\lambda_2,\lambda_3>0$, and $\lambda_{0}\geq 0$ such that $\lambda_1+\lambda_2+\lambda_3=1$.  
Hence, $U = U_0 +(\lambda_2\alpha,\lambda_3\beta)$ and by \cref{det:prop}, $\det(U-U_0)<0$. 
Finally, from the last equality in \cref{2conect} and \cref{quasidet},  we conclude that $U\notin Q^e(\mathcal{U})$ as claimed. 

\item Assume that there is only one compatible pair of wells in $\mathcal{U}$. 
By hypothesis, this pair is rank-one compatible, and without loss of generality let $U_3=U_1 +(a,0)$, and $U_2=U_1+(\xi,\eta)$, with $\xi\eta>0$ and $(\xi-a)\eta>0$.
It follows that $U_0:=U_1+(\xi,0)$ is rank-one compatible with all the wells in $\mathcal{U}$. 
From \cref{P:lamconv}, $L^e(\mathcal{U})=C(\{U_1,U_3\})\cup \{U_2\}$ and $C(\mathcal{U}) = L^e(\mathcal{U})\cup L$
where
\[
L = \left\{\lambda_1U_1 +\lambda_2U_2 +\lambda_3U_3\,|\, \lambda_1,\lambda_3\in[0,1],\ \lambda_2\in(0,1),\ \ \mbox{and}\ \ \lambda_1+\lambda_2+\lambda_3=1  \right\}.
\] 
Now, if $M\in L$ then $M=U_0 +(-\lambda_1\xi-\lambda_3(\xi-a),\lambda_2\eta)$. By a computation and \cref{det:prop} we get $\det(M-U_0)<0$. 
From the last inequality and \cref{quasidet}, we conclude that $M\notin Q^e(\mathcal{U})$, and the proof is complete.
\end{enumerate}
\end{proof} 

\subsection{The four-well case}

The strategy of the proof is to use the rank-one lines passing through some of ${\mathcal U}$'s well to split the affine space $\Pi_Q$. With this splitting, we consider all possible locations of the remaining wells. The  goal is to show that any matrix $M\in C(\mathcal{U})$ that belongs to the exterior of $L^e({\mathcal U})$ cannot be in $Q^e({\mathcal U})$. This painstaking procedure yields the result. 

\begin{proof}[Theorem~\ref{thm:4wells}]
Let $\Pi_Q$ be the affine space generated by $\mathcal{U}$, see \cref{planeQV}. 
Notice that, since $\{V_1,V_2\}$ and $\{W_1,W_2\}$ are rank-one compatible pairs, there are only two possibilities, either the corresponding rank-one lines generated by them are parallel or not.  
We divide the proof into three steps.

{\bf STEP 1} We assume  that $D=C(\{V_1,V_2\})\cup C(\{W_1,W_2\}$ is a disconnected set (\cref{it:le1} in \cref{thm:4wells}). The compatible planar cones $\mathscr{C}_Q(W_1)$ and $\mathscr{C}_Q(W_2)$ determine a division of the plane $\Pi_Q$ into three disjoint sets (see \cref{fig:nonconnec}), namely
\[
\Pi_Q=A \cup B\cup C, \quad
\]
where $A=\mathscr{C}_Q(W_1)\cap\mathscr{C}_Q(W_2)$, 
$B=\mathscr{C}_Q(W_1)\triangle\mathscr{C}_Q(W_2)$, and  $C=(\Pi_Q\setminus\mathscr{C}_Q(W_1))\cap(\Pi_Q\setminus\mathscr{C}_Q(W_2))$.
\begin{figure}[ht]
\begin{centering}
\includegraphics[scale=0.5]{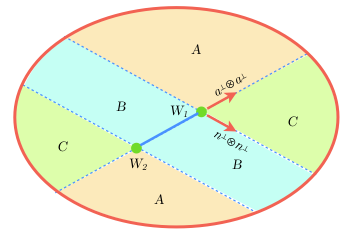}
\caption{The plane $\Pi_Q$ divided into three sets for step 1. This division is naturally determined by the rank-one lines through $W_1$ and $W_2$.}
\label{fig:nonconnec}
\end{centering}
\end{figure}
Notice that, $A$ is the set of all matrices compatible with $W_1$ and $W_2$ simultaneously, namely
\begin{equation}\label{eq:divA}
A = \left[\mathscr{C}_Q^-(W_1)\cap\mathscr{C}_Q^-(W_2)\right] \cup \left[\mathscr{C}_Q^+(W_1)\cap\mathscr{C}_Q^+(W_2)\right] \cup C(\{W_1,W_2\}).
\end{equation}
Since $D$ is disconnected,  $V_1$ and $V_2$ do not belong to $C(\{W_1,W_2\})$ or simultaneously to different $B$, or $C$, or $A$ components.
Next, we consider the remaining four possible cases.

Now, because $\{V_1,V_2\}$ is a rank-one compatible pair, if $V_1, V_2\in A$, both belong to $\mathscr{C}_Q^-(W_1)\cap\mathscr{C}_Q^-(W_2)$ or $\mathscr{C}_Q^+(W_1)\cap\mathscr{C}_Q^+(W_2)$. Therefore, all wells in $\mathcal{U}=\{W_1,W_2,V_1,V_2\}$ are pairwise compatible and $L^e(\mathcal{U}) =C(\mathcal{U})$. 
If one of the wells in $\{V_1,V_2\}$ belongs to $A$ and the other belongs to $B$, we  obtain the same conclusion.
 \medskip
 
 We assume the case when $V_1\in A$ and $V_2\in C$. 
 Thus,  $W_1,W_2$, and $V_2$ belong all to $\mathscr{C}_Q(V_1)$'s upper or lower parts by definition of $A$ (see \cref{eq:divA}). 
Without loss of generality, we assume that $W_1,W_2,V_2\in \mathscr{C}_Q^+(V_1)$, 
and let $W_2$ the furthest well in $\{W_1,W_2\}$ from $V_1$. 
By construction, $W_2= V_1 + (\xi,\eta)$, and $W_1= V_1 + (\xi',\eta)$ for some $\xi\geq 0 \geq \eta$, and $\xi \geq  \xi' \geq 0$, respectively. 
Now, we assume that $C(\{V_1,V_2\})$ and $C(\{W_1,W_2\})$ are non parallel segments, and 
$V_2=V_1+(0,\gamma)$ for some $\gamma\le 0$.
In this case, by computation, $U_0 = V_1 +(0,\eta)$ is such that $\det(U-U_0)\le 0$ for every $U\in \{W_1,W_2,V_1,V_2\}$.   
Now, let $M\in C(\mathcal{U})\setminus L^e(\mathcal{U})$. 
As  a consequence of \cref{thm:lamconv} and \cref{P:lamconv}, $M\in \relint C(\{W_2,U_0,V_2\})\cup C(\{W_2,V_2\})$.
Thus, $M= U_0 + (\alpha', \beta')$ for some $\beta'> 0 > \alpha'$ and $\det(M-U_0)<0$. 
From \cref{quasidet} we conclude that $M\notin Q^e(\mathcal{U})$. 
If $C(\{V_1,V_2\})$ and $C(\{W_1,W_2\})$ are parallel segments, 
letting $U_0 = V_1 +(\xi,0)$, we can repeat the argument to get the result. 

\medskip

Let $V_1$ and $V_2$ belong to $B$. 
 By definition of $B$, $V_1$ and $V_2$  are compatible with only one of the wells $W_1$ or $W_2$. 
The result follows by replacing wells $\{V_1,V_2\}$ with $\{W_1,W_2\}$ and arguing as in the previous paragraph.

\medskip
 
Next, we assume that $V_1\in B$ and $V_2\in C$. 
By definition of $B$, only one wells in $\{W_1,W_2\}$ is compatible with $V_1$, say $W_1$. Then, we have that
 \begin{equation}
 \begin{split}\label{eq:condBC}
 \det(W_1-W_2) =0,\quad \det(W_2-V_1)>0, \quad \det(V_1-W_1)\leq 0, \\
  \det(V_1-V_2) =0,\quad \det(V_2-W_1)>0, \quad \det(W_1-V_1)\leq 0.
 \end{split}
 \end{equation}
Thus, by \cref{thm:lamconv},  $L^e(\mathcal{U}) = L^e(\{V_2,V_1,W_1\}) \cup L^e (\{V_1,W_1,W_2\})$, and  by \cref{P:lamconv}, 
 \[
 L^e(\{V_2,V_1,W_1\}) = C(\{U_0, V_1,W_1\}) \cup C(\{U_0, V_2\}),
 \] 
 and
 \[ L^e(\{V_1,W_1,W_2\}) = C(\{U_0', V_1,W_1\}) \cup C(\{U_0', W_2\}),
 \]
 where $U_0\in C(\{V_1, V_2\})$ and $U_0'\in  C(\{W_1, W_2\})$ are such that $\det(W_1-U_0) = 0$ and $\det(V_1-U_0') = 0$, respectively. Hence, 
$
C(\mathcal{U}) = L^e(\mathcal{U})\cup L_1 \cup L_2
$ 
where $L_1 =\relint C(\{V_1,V_2,W_1\}) \cup C(\{V_2,W_1\})$ and $L_2= \relint C(\{V_1,W_1,W_2\}) \cup C(\{V_1,W_2\})$.
By construction, every $M\in L_1$ and $N\in L_2$ satisfy that $\det(M-U_0)<0$ and $\det(N-U_0')<0$. From last inequalities, \cref{eq:condBC}, and \cref{quasidet}, we conclude $M,N\not\in Q^e(\mathcal{U})$ as claimed.
 
\medskip

\begin{figure}[ht]
     \centering
     \begin{subfigure}[b]{0.3\textwidth}
         \centering
         \includegraphics[width=\textwidth]{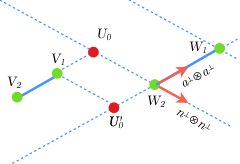}
         \caption{\footnotesize $C(\{V_1, V_2\}) \parallel C(\{W_1, W_2\}) $}
         \label{fig:4wells-CP}
     \end{subfigure}
     \hspace{2cm}
     \begin{subfigure}[b]{0.3\textwidth}
         \centering
         \includegraphics[width=\textwidth]{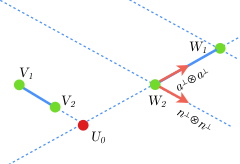}
         \caption{\footnotesize $C(\{V_1, V_2\}) \nparallel C(\{W_1, W_2\})$}
         \label{fig:4wells-CO}
     \end{subfigure}
     \caption{Figures (a) and (b) are two four-well configurations considered in step 1 where two wells are incompatible with the remaining pair}
     \label{fig:Twographs}
\end{figure}

Assume  $V_1, V_2\in C$.
Since $D$ is disconnected,  by \cref{thm:lamconv}, $L^e(\mathcal{U}) = D$. First, we assume that $C(\{V_1,V_2\})$ and $C(\{W_1,W_2\})$ are non parallel segments, see \cref{fig:4wells-CO}. 
Let $U_0$ be the intersection point of  the rank-one lines $\{tV_1 + (1-t)V_2\,|\, t\in \re\}$ and $\{tW_1 + (1-t)W_2\,|\, t\in \re\}$. 
We assume that $W_1-U_0=(\alpha,0)$ and $V_1-U_0=(0,\beta)$ for some $\alpha,\, \beta\in \re$ such that $\alpha\beta<0$, since $V_1$ and $W_1$ are incompatible.    
In this case 
$C(\mathcal{U}) \subset \{U_0 + t(W_1-U_0) +s(V_1-U_0)\,|\, t,s\in \re^+\cup\{0\} \}$, see \cref{fig:4wells-CO}. As before, we conclude that $\det(U-U_0)=0$ for every $U\in \mathcal{U}$.  Hence, if $M\in C(\mathcal{U})\setminus D$, $\det(M-U_0)<0$ and by \cref{quasidet}, $M\notin Q^e(\mathcal{U})$.

Second, we assume that $C(\{V_1,V_2\})$ and $C(\{W_1,W_2\})$ are parallel segments, 
see \cref{fig:4wells-CP}. 
Without loss of generality, let $V_2=V_1+(\xi,0)$ for some $\xi\neq 0$ and let $V_1$, $W_2$ be the closest wells between $\{V_1,V_2\}$ and $\{W_1,W_2\}$, then $W_2=V_1+(\alpha,\beta)$ where
$\alpha\beta > 0$ and $(\alpha-\xi)\beta >0$  since $V_1\,,V_2\in C$. Moreover, $W_1=V_1+(\alpha+\gamma,\beta)$ where $\gamma\alpha>0$ since $W_2$ is closer to $V_1$ than $W_1$.

Now, let $U_0=V_1+(\alpha,0)$ and $U_0'=V_1+(0,\beta)$.  A calculation yields, 
\begin{equation}\label{eq:conectCC}
\det(U-U_0)\ge 0,\quad\text{ and }\det(U-U_0')\ge 0
\end{equation}
for every $U\in \mathcal{U}$. 
If $M\in  \relint  \mathscr{C}_Q^+(U_0') \setminus D$, $\det(M-U_0')>0$.
Thus, last inequality, \eqref{eq:conectCC}, and 
\cref{quasidet} yields $M\notin Q^e(\mathcal{U})$.
Analogously, if  $M\in  \relint  \mathscr{C}_Q^-(U_0) \setminus D$ then $\det(M-U_0)>0$ and we conclude that $M\notin Q^e(\mathcal{U})$.
Finally, we notice that
\[
C(\{V_1,W_1,W_2\})\subset \mathscr{C}_Q^+(U_0')
\quad\text{ and }\quad
C(\{W_2,V_1,V_2\})\subset \mathscr{C}_Q^-(U_0),
\] 
and $C(\mathcal{U})=C(\{V_1,W_1,W_2\})\cup C(\{W_2,V_1,V_2\})$. Hence, if $M\in C(\mathcal{U})\setminus D$, then $M\notin Q^e(\mathcal{U})$ and the result follows.

{\bf STEP 2} 
We assume that the intersection of the sets $\{V_1,V_2\}$ and $\{W_1,W_2\}$ has only one element, say $V$, and $D$ is contained either in the upper or in the lower part of $\mathscr{C}_Q(V)$
 (\cref{it:le2} in \cref{thm:4wells}). Let $\mathcal{U}=\{U_1,U_2,U_3,U_3\}$, assume $V=U_2$, $U_4\notin \{V_1,V_2\}\cup\{W_1,W_2\}$, and $D\subset \mathscr{C}_Q^+(U_2)$ (see \cref{fig:4wells-V}). Thus, $U_1=U_2+(\alpha,0)$, $U_3=U_2+(0,\beta)$, and $U_4 = U_2 +(l,m)$ for some $\beta<0<\alpha$ and $l,\, m\in \re$.
 Notice that  $U_4\notin \mathscr{C}^-_Q(U_2)$, otherwise  $l<0<m$ and  $\mathcal{U}$ would be a wedge configuration. 
 Now,  we split $\Pi_Q\setminus \mathscr{C}^-_Q(U_2)$ into five regions by the rank-one lines passing through  $U_1,\,U_2$ and $U_3$, see \cref{fig:4wells-V}. Namely, $\Pi_Q\setminus \mathscr{C}^-_Q(U_2)=A\cup B\cup C\cup D \cup E$, where 
 \[
\begin{array}{c} A =\{U\in\Pi_Q\setminus \mathscr{C}^-_Q(U_2)\,|\, U\in \mathscr{C}^+_Q(U_2) \cap \mathscr{C}_Q(U_1) \cap \mathscr{C}_Q(U_3)\}\\
B =\{U\in\Pi_Q\setminus \mathscr{C}^-_Q(U_2)\,|\, U\in \mathscr{C}^+_Q(U_2) \cap \mathscr{C}_Q(U_1) \triangle \mathscr{C}_Q(U_3)\}\\
C =\{U\in\Pi_Q\setminus \mathscr{C}^-_Q(U_2)\,|\, U\in \mathscr{C}^c_Q(U_2) \cap \mathscr{C}_Q(U_1) \triangle \mathscr{C}_Q(U_3)\}\\
D =\{U\in\Pi_Q\setminus \mathscr{C}^-_Q(U_2)\,|\, U\in \mathscr{C}^+_Q(U_2) \cap \mathscr{C}^c_Q(U_1) \cap \mathscr{C}^c_Q(U_3)\}\\
E =\{U\in\Pi_Q\setminus \mathscr{C}^-_Q(U_2)\,|\, U\in \mathscr{C}^c_Q(U_2) \cap \mathscr{C}^c_Q(U_1) \cap \mathscr{C}^c_Q(U_3)\}
\end{array} 
\]  
and every complement set operation considers $\Pi_Q\setminus \mathscr{C}^-_Q(U_2)$ as the universe. 
\begin{figure}[ht]
\begin{centering}
\includegraphics[scale=0.5]{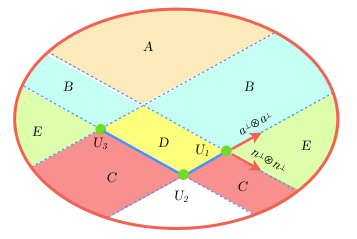}
\caption{This figure shows the splitting of $\Pi_Q$ used in step 2.}
\label{fig:4wells-V}
\end{centering}
\end{figure}

Next, we explore $U_4$'s five different possibilities. If $U_4\in A$, all wells are pairwise compatible and $L^e(\mathcal{U})=Q^e(\mathcal{U})=C(\mathcal{U})$ (see \cref{rel-boun-L2}). 
By hypothesis, $U_4\notin E$ otherwise it will be incompatible with the remaining wells in $\mathcal{U}$. 

Assume $U_4\in B$, thus it is incompatible with either $U_1$ or $U_3$. Without loss of generality, let $\det (U_4-U_3)>0$, that is  $l(m-\beta)>0$.  
Hence,  $L^e(\mathcal{U}) = L^e(\{U_4,U_2,U_3\}) \cup L^e(\{U_4,U_2,U_1\})$ by \cref{thm:lamconv}. Also by \cref{P:lamconv}, we accomplish that $L^e(\{U_4,U_2,U_1\}) = C(\{U_1,U_2,U_4\})$, and  $L^e(\{U_4,U_2,U_3\}) = C(\{U_2,U_3\})\cup C(\{U_4,U_2,U_0\})$. 
Here,  $U_0 = U_2 +(0,m)$ for some $\beta\leq m\leq 0<l$ is the flag's center of  $\{U_2,U_3,U_4\}$. Thus, $\det(U-U_0)=0$ for $U\in\{U_2,U_3,U_4\}$, and $\det(U_1-U_0) =-\alpha m|a\times n|^2\ge 0$ (see \cref{det:prop}).

Now, since $C(\mathcal{U}) = C(\{U_1,U_2,U_4\})\cup C(\{U_2,U_3,U_4\})$ it follows that $C(\mathcal{U}) = L^e(\mathcal{U}) \cup L$, where 
$$L = \{\lambda_1U_0+\lambda_2U_4+\lambda_3U_3\,|\, \lambda_1\in[0,1),\ \lambda_2,\lambda_3\in(0,1), \mbox{ and } \lambda_1+\lambda_2+\lambda_3=1\}.$$ 
If $U\in L$,  $U=U_0+(\lambda_2l,\lambda_3(\beta-m))$ for $ \lambda_2,\lambda_3\in(0,1)$, and  $\det(U-U_0)<0$ by \cref{det:prop}. Therefore, by \cref{quasidet}, $U\notin Q^e(\mathcal{U})$ as claimed. 

Let $U_4\in C$. Then, $\det(U_2-U_4)>0$ and without loss of generality, we assume $\det(U_4-U_3)>0\geq \det(U_4-U_1)$. 
By \cref{thm:lamconv}, we get $L^e(U) = L^e(\{U_1,U_2,U_3\})\cup L^e(\{U_1,U_2,U_4\})$.
Due to \cref{P:lamconv} and the compatibility relations among $\mathcal{U}$'s wells, we have that $L^e(\{U_1,U_2,U_3\})=C(\{U_1,U_2\})\cup C(\{U_2,U_3\})$  and $L^e(\{U_1,U_2,U_4\}) = C(\{U_0,U_1,U_4\})\cup C(\{U_0,U_2\})$, where $U_0\in  C(\{U_1,U_2\})$ is the flag's center of $\{U_1,U_2,U_4\}$. 
Hence, $C(\mathcal{U}) = L^e(\mathcal{U})\cup  L_1\cup L_2$, 
where 
\[
\begin{split}
L_1=\left\{\lambda_1U_1+\lambda_2U_2+\lambda_3U_3 \,|\, \lambda_1,\lambda_3\in (0,1), \lambda_2\in[0,1) \ \mbox{and} \ \lambda_1+\lambda_2+\lambda_3=1\right\},\\
L_2=\left\{\lambda_1U_0+\lambda_2U_2+\lambda_3U_4 \,|\, \lambda_2,\lambda_3\in (0,1), \lambda_1\in[0,1) \ \mbox{and} \ \lambda_1+\lambda_2+\lambda_3=1\right\}.
\end{split}
\]
If $M\in L_1$, then $M=U_2 +(\lambda_1\alpha,\lambda_3\beta)$, and
$\det(U-U_2)=\lambda_1\lambda_2\alpha\beta|a\times n|^2 < 0$. Thus,
recalling that $\det(U-U_2)\ge 0$ for every $U\in \mathcal{U}$, we get $M\notin Q^e(\mathcal{U})$ by \cref{quasidet}. Now, if $M\in L_2$, we can prove by the same arguments that $\det(M-U_0)<0$. 
\medskip

Finally, let $U_4\in D$. In this case, $U_4$ is compatible only with $U_2$, and,  by \cref{thm:lamconv}, $L^e(\mathcal{U}) = L^e(\{U_2,U_3,U_4\})\cup L^e(\{U_1,U_2,U_4\})$, and by \cref{P:lamconv}, 
\[
\begin{array}{l}
L^e(\{U_2,U_3,U_4\}) = C(\{U_0,U_2,U_4\})\cup C(\{U_0,U_3\}), \text{ and }\\
L^e(\{U_1,U_2,U_4\}) = C(\{U'_0,U_2,U_4\})\cup C(\{U'_0,U_1\}),  
\end{array}
\]
where $U_0\in C(\{U_2,U_3\})$, $U'_0\in C(\{U_1,U_2\})$, and $\det(U_4-U_0) = \det(U_4-U'_0)=0$. It follows that $U_0 = U_2 +(0,m)$ and $U'_0 = U_2 +(l,0)$ for some $l,m\in \re$ such that $lm<0$. 
Due to the location of $U_4$,  $C(\mathcal{U})$ is either a quadrilateral or a triangle. In the former case, $C(\mathcal{U})=C(\{U_1,U_2,U_4\})\cup C(\{U_2,U_3,U_4\})$, so  $C(\mathcal{U}) = L^e(\mathcal{U}) \cup L_1 \cup L_2$, where
\[
L_1=\left\{\theta_0U_0+\theta_3U_3 +\theta_4U_4 \,|\, \theta_3,\theta_4\in (0,1), \theta_0\in[0,1) \ \mbox{and} \ \theta_0+\theta_3+\theta_4=1\right\},
\]
\[
L_2=\left\{\lambda_1U_1+\lambda_0U'_0 +\lambda_4U_4 \,|\, \lambda_1,\lambda_4\in (0,1), \lambda_0\in[0,1) \ \mbox{and} \ \lambda_0+\lambda_1+\lambda_4=1\right\}.
\]
If $M\in L_1$, then $\det(M-U'_0)<0$ since $M=U'_0 +(\lambda_1(\alpha-l),\lambda_4m)$ and $m(\alpha-l)<0$, due to 
 the incompatibility between $U_4$ and $U_1$.  Thus, by \cref{quasidet}, $L_1\cap Q^e(\mathcal{U})=\emptyset$. The same result follows  for $U\in L_2$ by a similar  argument. 

Now we assume that $C(U)$ is a triangle, hence $U_4\in C(\{U_1,U_2,U_3\})$. 
As before, we let $U_0 = U_2 +(0,m)$ and $U'_0 = U_2 +(l,0)$, and we define $U'_1,\,U'_3\in C(\{U_1,U_3\})$ such that $U'_1 = U_2+(t,m)$ and $U'_3 = U_2 +(l,s)$ for some $t,s \in \re$, see \cref{fig:caseD}. By construction, it follows that $C(\mathcal{U})\subset L^e(\mathcal{U})\cup \tilde{L}_1\cup \tilde{L}_2$, where 
\[
\tilde{L}_1=\left\{\theta_0U_0+\theta_1U'_1 +\theta_3U_3  \,|\, \theta_1,\theta_3\in (0,1), \theta_0\in[0,1) \ \mbox{and} \ \theta_0+\theta_1+\theta_3=1\right\},
\]
\[
\tilde{L}_2=\left\{\lambda_0U'_0 +\lambda_1U_1+\lambda_3U'_3 \,|\, \lambda_1,\lambda_3\in (0,1), \lambda_0\in[0,1) \ \mbox{and} \ \lambda_0+\lambda_1+\lambda_3=1\right\}.
\]
The result follows as in the quadrilateral case by replacing $U_4$ with $U'_1$ and $U_4$ by $U'_3$ in the sets $L_1$, and $L_2$, respectively.

\begin{figure}[ht]
\begin{centering}
\includegraphics[scale=0.5]{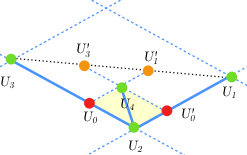}
\caption{This figure shows a four-well configuration considered in step 2, when $U_4\in D$ and $C(\mathcal{U})$ is a triangle.}
\label{fig:caseD}
\end{centering}
\end{figure}

{\bf STEP3}: We assume $D$ is a connected set and ${\mathcal U}\subset D$ ( \cref{it:le3} in \cref{thm:4wells}). In this case, 
the rank-one segments $C(\{V_1,V_2\})$ and $C(\{W_1,W_2\})$ intersects at a point, say $U_0$. Without loss of generality, assume that $U_0 = W_1 +(-\alpha,0)$ for some $\alpha>0$. Hence, by \cref{planeQV}, we get
\[
W_1 = U_0+(\alpha,0), \quad W_2 = U_0+(\beta,0), \quad V_1 = U_0 +(0,\gamma), \quad \mbox{and} \ \ V_2 = U_0 +(0,\eta)  
\]
where $\alpha \beta \leq0$ and $\gamma \eta\leq 0$. Due to  \cref{det:prop} we conclude that $W_1$ is compatible with either $V_1$ or $V_2$. Assume $\det(W_1-V_1)\leq 0$ and $\det(W_2-V_2)\leq 0$, then $L^e(\mathcal{U}) = L^e(\{W_1,W_2,V_1\})\cup L^e(\{W_1,W_2,V_2\})$ by  \cref{thm:lamconv}. Moreover, from \cref{P:lamconv}, $L^e(\{W_1,W_2,V_1\}) = C(\{W_1,V_1,U_0\})\cup C(\{W_2,U_0\})$ and $
L^e(\{W_1,W_2,V_2\}) = C(\{W_2,V_2,U_0\})\cup C(\{W_1,U_0\})$. Thus, $C(\mathcal{U}) = L^e(\mathcal{U})\cup L_1\cup L_2$ where 
\[
L_1=\left\{\theta_0U_0+\theta_1V_1 +\theta_3W_2  \,|\, \theta_1,\theta_3\in (0,1), \theta_0\in[0,1) \ \mbox{and} \ \theta_0+\theta_1+\theta_3=1\right\}, \text{ and }
\]
\[
L_2=\left\{\lambda_0U_0+\lambda_1W_1 +\lambda_3V_2  \,|\, \lambda_1,\lambda_3\in (0,1), \lambda_0\in[0,1) \ \mbox{and} \ \lambda_0+\lambda_1+\lambda_3=1\right\}.
\]
As before, a calculation yields $\det(U-U_0)<0$ for every $U\in L_1\cup L_2$,  but $\det(V-U_0)=0$ for all $V\in \mathcal{U}$. Therefore, $(L_1\cup L_2)\cap Q^e(\mathcal{U}) = \emptyset$ and the result follows. 
\medskip
 
As we have considered all admissible configurations, this concludes the theorem's proof.
\end{proof}

\subsection{The \texorpdfstring{$n$}{}-well problem (basic configurations)}

In this final section, we prove \cref{thm:iterative}.
We begin with  two preliminary lemmas.

\begin{lemma}\label{prop:bc}
A finite coplanar set $\mathcal{U}$ is a basic configuration made of $n$ basic blocks if and only if there exist $M_0\in \Pi_Q$, and $\alpha_i, \beta_i$ for $i=1,2,\cdots n+1$ positive constants such that if 
\begin{equation}\label{Vis}
V_i=M_0 +\sum_{j<i}\left(\alpha_j,\beta_j\right)+(\alpha_i,0),\quad \mbox{and}\quad  W_i=M_0 +\sum_{j<i}\left(\alpha_j,\beta_j\right)+(0,\beta_i),
\end{equation}
then 
$\left(\bigcup_{i=1}^{n+1}\{V_i,W_i\}\right)\setminus \mathcal{U} \subset \{A,B\}$,  where 
$A\in \{V_1,W_1\}$ and $B\in \{V_{n+1},W_{n+1}\}$. 
\end{lemma}
\begin{remark}
We emphasize that $\bigcup_{i=1}^{n+1}\{V_i,W_i\}$ can be equal to $\mathcal{U}$, or their set difference has at most two wells. 
\end{remark}
\begin{proof}
We assume $\mathcal{U}$ is a basic configuration to prove the result by induction on the number of basic blocks in $\mathcal{U}$.
The induction's base is for $n=1$. Let $\mathcal{U}$ be a basic block with center at $U_0=(0,0)$ (see \cref{fig:bb3} and \cref{fig:bb4}). 
First, we assume that  $\mathcal{U}$ is a four-well basic block. 
Since all wells in $\mathcal{U}$ are rank-one compatible with $U_0$, by \cref{planeQV},  $\mathcal{U}=\{(\alpha,0),(0,\beta), (\alpha',0),(0,\beta')\}$ for some $\alpha>0>\alpha'$ and  $\beta>0>\beta'$. 
By letting 
\[
M_0=(\alpha',\beta'), \quad V_1 = (0,\beta'), \quad W_1 = (\alpha',0), \quad V_2 = (\alpha,0),\text{ and } W_2 = (0,\beta),
\]
we have that $\mathcal{U}=\{V_1,W_1,V_2,W_2\}$, and \cref{Vis} holds for the positive constants $\alpha_1 = -\alpha'$, $\beta_1 = -\beta'$, $\alpha_2 = \alpha$ and $\beta_2 = \beta$.  
Second, assume $\mathcal{U}$ is a three-well basic block. Since all wells in $\mathcal{U}$ are rank one compatible with the flag's center $U_0=(0,0)$, we have that $\mathcal{U}=\{(\alpha,0),(0,\beta),\alpha'(\chi,1-\chi) \}$, for some  $\chi\in\{0,1\}$ and, $\alpha\beta>0>\alpha\alpha'$. 
Next, if $\beta'$ is such that $\beta\beta'<0$, then the set $\mathcal{U}'= \mathcal{U}\cup\{\beta'(1-\chi,\chi)\}$ is a four-well basic block. Now, the result follows by the above argument for the four-well case.

For the induction step, we assume that the result holds for every basic configuration made of $n$ basic blocks and we will prove that it is still valid for a basic configuration $\mathcal{U}$ with $n+1$ basic blocks. 
Indeed $\mathcal{U} = \mathcal{U}'\cup \mathcal{U}_{n+1}$, where $\mathcal{U}'$ is a $n$ basic block configuration. 
Thus, by the induction hypothesis, $\bigcup_{i=1}^{n+1}\{V_i,W_i\} \setminus \mathcal{U}'\subset \{A,B\}$. 
By construction, the set $\mathcal{U}_n$ has two adjacent basic blocks, ($\mathcal{U}_{n-1}$ and $\mathcal{U}_{n+1}$), and  it has two different compatible, but not rank-one compatible, pairs of wells.
Since, only four-well basic block has this property, we have that $\bigcup_{i=1}^{n+1}\{V_i,W_i\} \setminus \mathcal{U}'\subset \{A\}$. 

Next, since no more than two wells in $\mathcal{U}$ are colinear, $\mathcal{U}_{n+1}$ is either a four-well or three-well basic block. In the former case, the set $\mathcal{U}_{n+1}=\{V_{n+1},W_{n+1},V_{n+2},W_{n+2}\}$ where  $V_{n+1}$, $W_{n+1}$ satisfies \cref{Vis} and 
\[
 V_{n+2} = \sum_{i=1}^{n+1}\left(\alpha_i,\beta_i\right) +(\alpha_{n+2},0), \quad \mbox{and} \quad  W_{n+2} = \sum_{i=1}^{n+1}\left(\alpha_i, \beta_i\right) +(0,\beta_{n+2}) 
\]
for some positive $\alpha_{n+2}$ and $\beta_{n+2}$ as claimed. If $\mathcal{U}_{n+1}$ is a three-well basic block, then $\mathcal{U}_{n+1}=M_0^{n+1}+\{ (-\alpha_{n+1},0), (0,-\beta_{n+1}), \alpha'(\chi,1-\chi) \}$ with $\alpha'>0$ and $M_0^{n+1}=\sum_{i=1}^{n+1}(\alpha_i,\beta_i)$. 
As in the proof of the induction basis, we complete the set $\mathcal{U}_{n+1}$ to form a four-well basic block, and we conclude by the four-well case. 

\medskip

Notice that if we prove that $\mathcal{U}_{i}=\{V_i,W_i,V_{i+1},W_{i+1}\}$ is a basic block for $i=1,2,3,\cdots, n$, the reverse implication, namely $ \bigcup_{i=1}^{n+1}\{V_i,W_i\}\setminus \mathcal{A}$ is a basic configuration for each set $\mathcal{A}\subset \{A,B\}$, follows straight forward since  $\mathcal{U}_1\setminus \{A\}$ and $\mathcal{U}_n\setminus \{B\}$ also are three-well basic blocks. 

Hence, we prove that $\mathcal{U}_{i}$ is a basic block for $i=1,2,3,\cdots, n$. Indeed, since $W_{i+1}-V_i=(0,\beta_i+\beta_{i+1})$  and  $V_{i+1}-W_i=(\alpha_i+\alpha_{i+1},0)$, see \cref{Vis}, we have that $\{V_i, W_{i+1}\}$  and $\{V_{i+1}, W_{i}\}$ are both rank-one compatible pairs. Moreover,  the well $M_0 +\sum_{j\leq i}\left(\alpha_i,\beta_i\right)$ belongs to both $C(\{V_i, W_{i+1}\})$ and $C(\{V_{i+1}, W_{i}\})$.  Hence, $\mathcal{U}_i$ satisfies the conditions of \cref{it:le3} \cref{thm:4wells}, and it is a four-well basic block and the proof is complete. 
\end{proof}

\begin{lemma}\label{lem:supplines}
Let $\mathcal{U}$ be a basic configuration such that $\mathcal{U}=\bigcup_{i=1}^{n}\mathcal{U}_i$ where $\mathcal{U}_i$ is a basic block for every $i=1,2,\cdots,n$. Then there exist two compatible wells $P_1,Q_1\in \mathcal{U}_1$, and two compatible wells $P_n,Q_n\in \mathcal{U}_{n}$ such that $\Aff(\{P_1,Q_1\})$ and $\Aff(\{P_n,Q_n\})$ are two supporting lines of $C(\mathcal{U})$
\end{lemma}
\begin{proof}
To keep the proof simple, let $M_0=0$. 
We know that $\mathcal{A}=\bigcup_{i=1}^{n+1}\{V_i,W_i\}\setminus \mathcal{U}$ is empty or it has at most two elements. We assume either $\mathcal{A}=\{A\}$ or $\mathcal{A}=\emptyset$. In the latter case, we let $(P_1,Q_1)=(V_1,W_1)$, and define the linear functional $\ell(U):=\nu_1\alpha+\nu_2\beta$, 
where $\nu_1=1/\alpha_1$,   $\nu_2=1/\beta_1$, and $U=(\alpha,\beta)$. Thus, 
\[
C(\{P_1,Q_1\}) \subset \Aff(\{P_1,Q_1\}) =\{U\in \Pi_Q \,| \,\ell(U)= 1\}.
\]
Next, since every $\alpha_i,\beta_i$ in \eqref{Vis} are positive numbers, 
\begin{equation}\label{eq:vbounds}
1=\ell(V_1)\leq\ell(V_{i}),\ \mbox{and} \ 1=\ell(W_1)\leq \ell(W_{i}),
\end{equation}
for $i=1,2,\cdots n+1$. Thus, all wells in $\mathcal{U}$ are on one side of the line $\ell(U)=1$, and $\Aff(\{P_1,Q_1\})$  is a supporting line of  $C(\mathcal{U})$. Assuming $\mathcal{A}=\{A\}$, we have two further cases, either $A=V_1$ or $A=W_1$. We let $(P_1,Q_1)=(V_2,W_1)$ with $\nu_1=0$,    $\nu_2=1/\beta_1$ and $(P_1,Q_1)=(V_1,W_2)$ with $\nu_1=1/\alpha_1$, $\nu_2=0$, respectively. Now the proof follows exactly as above.

In the case of $\mathcal{A}=\emptyset$, we also let $(P_n,Q_n)=(V_{n+1},W_{n+1})$ and $\tilde{\ell}(U)=\tilde{\nu}_1 \alpha+ \tilde{\nu_2}\beta$ with 
$
\tilde{\nu}_1= \beta_{n+1}/\kappa$   and 
$ \tilde{\nu}_2= \alpha_{n+1}/\kappa$,
where $\kappa = \alpha_{n+1}\beta_{n+1} +\sum_{i=1}^{n}\alpha_{n+1}\beta_i+\beta_{n+1}\alpha_i$. Thus, 
\[
C(\{P_n,Q_n\}) \subset \Aff(\{P_n,Q_n\}) =\{U\in \Pi_Q \,| \,\tilde{\ell}(U)= 1\},
\]
and 
\begin{equation}\label{eq:vbounds2}
\tilde{\ell}(V_i)\leq\tilde{\ell}(V_{n})=1,\ \mbox{and} \ \tilde{\ell}(W_i)\leq \tilde{\ell}(W_{n})=1,
\end{equation}
for $i=1,2,\cdots n+1$. Arguing as before, we find that $\Aff(\{P_n,Q_n\})$ is a supporting line of $C(\mathcal{U})$.
Assuming $\mathcal{A}=\{B\}$, we have two cases, either $B=V_{n+1}$ or $B=W_{n+1}$. 
If $B=V_{n+1}$, 
we let $(P_n,Q_n)=(V_n,W_{n+1})$ with  $\tilde\nu_1=1/\sum_{j\le n}\alpha_j$ and $\tilde\nu_2=0$. If $B=W_{n+1}$, then we let $(P_n,Q_n)=(V_{n+1},W_{n})$ with 
$\tilde\nu_1=0$ and $\tilde\nu_2=1/\sum_{j\le n}\beta_j$. In both cases the result follows as before. 
Finally, the case $\mathcal{A}=\{A,B\}$ follows as combination of above cases, and the proof is concluded.
\end{proof}

Now, we present the proof of \cref{thm:iterative}.

\begin{proof}[Theorem \ref{thm:iterative}]
Since $\mathcal{U}$ is a basic configuration, $\mathcal{U}\subset \bigcup_{i=1}^{n+1}\{V_i,W_i\}$ for $V_i$ and $W_i$ given by \cref{Vis} in \cref{prop:bc}.
Let $\alpha_i,\beta_i$ for $i=1,2,\cdots n+1$ be the corresponding positive constants, and let 
\begin{equation}\label{eq:U0n}
U_0^{i}:=M_0 +\sum_{j\leq i}(\alpha_j,\beta_j)\ \text{ for }\ i\in\{1,\dots,n\}.
\end{equation}

First, we claim that $\det(U-U_0^i)\geq 0$ for every $U\in\mathcal{U}$. Indeed, by \cref{Vis}, 
\[
U_0^{i}-V_k 
= \left(\sum_{j\leq i}\alpha_j-\sum_{j\leq k}\alpha_j,\sum_{j\leq i}\beta_j-\sum_{j<k}\beta_j\right)=
\begin{cases}
\sum_{k<j\leq i} \left(\alpha_j, \beta_j\right) & \mbox{if } k<j,\\
\left(0,\beta_k\right) & \mbox{if } k=j,\\
-\sum_{i<j\leq k} \left(\alpha_j,\beta_j\right) & \mbox{if } k<j,
\end{cases}
\]
and 
\[
U_0^{i}-W_k = \left(\sum_{j\leq i}\alpha_j-\sum_{j < k}\alpha_j,\sum_{j\leq i}\beta_j-\sum_{j\leq k}\beta_j\right)
=\begin{cases}
\sum_{k\leq j\leq i}\left(\alpha_j,\beta_j\right) & \mbox{if } k<j,\\
\left(\alpha_k,0\right) & \mbox{if } k=j,\\
-\sum_{i<j< k}\left(\alpha_j,\beta_j\right) & \mbox{if } k<j.
\end{cases}
\]
Hence, $\det (U_0^{i}-V_k)\geq 0$ and $\det (U_0^{i}-W_k)\geq 0$ for every $i\in\{1,2,\cdots n\}$ and $k\in\{1,2,\cdots n+1\}$ by \cref{det:prop}, and the claim follows.
\medskip



Second, we show next that if $U\in C(\mathcal{U})$  and it does not belong to the lamination convex hull of any basic block, then $U\notin Q^e(\mathcal{U})$.
If $\mathcal{U}_i$ for $i=0,\dots, n$ are the  basic blocks of $\mathcal{U}=\bigcup_{i=1}^n\mathcal{U}_i$,
then, 
\begin{equation}\label{bcontentions}
\bigcup_{i=1}^n L^e(\mathcal{U}_i)\subset L^e(\mathcal{U})\subset Q^e(\mathcal{U})\subset C(\mathcal{U})\subset \Pi_Q. 
\end{equation}
Now, we assume that $U\in C(\mathcal{U})\setminus \bigcup_{i=1}^n L^e(\mathcal{U}_i)= \bigcap_{i=1}^n \left(C(\mathcal{U})\setminus L^{e}(\mathcal{U}_i)\right)$. For the sets in the latter intersection we have 
\[C(\mathcal{U})\setminus L^{e}(\mathcal{U}_i)=\mathcal{A}_i \cup \mathcal{B}_i, \ \ \mbox{with } \ \
\begin{array}{l}
\mathcal{A}_i=\{U\in C(\mathcal{U})\,|\, \det(U-U_0^i)<0 \}, \\ \mathcal{B}_i=\{U\in C(\mathcal{U})\,|\, \det(U-U_0^i)\geq 0,\ \ U\notin L^{e}(\mathcal{U}_i)\},
\end{array}
\]
where $U_0^i\in C(\mathcal{U}_i)$ (see \eqref{eq:U0n}) is the center of the basic block $\mathcal{U}_i$. We claim that $\bigcap_{i=1}^n\mathcal{B}_i=\emptyset$. 
Assuming the claim, by set algebra
\[\bigcap_{i=1}^n \left(C(\mathcal{U})\setminus L^{e}(\mathcal{U}_i)\right) = \bigcup_{i=1}^n\left(\bigcap_{j=1}^i\mathcal{A}_j\bigcap_{k>i}^n\mathcal{B}_k\right).
\]
Thus, there exists $i\in \{1,2,\cdots,n\}$ such that $\det(U-U_0^i)<0$. 
Hence, by the first claim and \cref{quasidet}, the result follows and the proof is complete.

Finally, we prove the claim. By \cref{lem:supplines} there exist two compatible wells $P_1,Q_1\in \mathcal{U}_1$, and another pair of compatible wells $P_n,Q_n\in \mathcal{U}_n$ such that $\Aff\{P_1,Q_1\}$ and $\Aff\{P_n,Q_n\}$ are two supporting lines of $C(\mathcal{U})\supset \cup_{i=1}^nL^e(\mathcal{U}_i)$. Equivalently, there exist two linear functionals $\ell, \tilde{\ell}:\Pi_Q\rightarrow \re$, see \cref{eq:vbounds} and \cref{eq:vbounds2}, such that $\ell(P_1) \leq \ell(U)$ and $\tilde{\ell}(U) \leq \tilde{\ell}(P_n)$  for every $U\in C(\mathcal{U})$. 

By contradiction, if  $U\in \bigcap_{i=1}^n\mathcal{B}_i\neq \emptyset$, 
then $U = U_0^i + (\gamma_i, \delta_i)$ for some $\gamma_i\delta_i\geq 0$, and $U\notin L^e(\mathcal{U}_i)$ for every $i=1,2,\cdots,n$. Now, we have two options: (a) $U\notin L^e(\mathcal{U}_n)$ and $U = U_0^n + (\gamma_n, \delta_n)$ for some $\gamma_n\geq 0$ and $\delta_n\geq 0$ or (b) $U\notin L^e(\mathcal{U}_1)$ and $U = U_0^1 + (\gamma_1, \delta_1)$ for some $\gamma_1\leq 0$ and $\delta_1\leq 0$. Since  $\mathcal{A}_1\cap \mathcal{B}_1=\emptyset$, $\mathcal{A}_n\cap \mathcal{B}_n=\emptyset,$ and $\mathcal{U}$ is a basic configuration, we conclude 
that $\ell(U)<\ell(P_1)$ and $\tilde{\ell}(U)>\tilde{\ell}(P_n)$ in cases (a) and (b), respectively. So
we get a contradict to $U\in C(\mathcal{U})$ and the claim is proved.



\end{proof}


\section{Acknowledgements}
AC was partially founded by CONACYT CB-2016-01-284451 and COVID19 312772 grants and a RDCOMM and grant UNAM PAPPIT–IN106118 grant. LM was partially founded by UNAM PAPPIT–IN106118 grant. Also, LM acknowledge the support of CONACyT grant--590176 during his graduate studies.

%




\bibliographystyle{abbrv} 
\cleardoublepage
\bibliography{MyRef2} 
\end{document}